\def\WHO{nobody} 
\def\version{30.09.2022}
\def\users{us}  %
\def\users{final-layout}   
\documentclass[12pt,leqno]{article}
\usepackage{upgreek}
\usepackage{xcolor}
\usepackage{bm,amsmath,amsthm,hyperref,amsfonts,amssymb,color}
\usepackage{mathrsfs} 
\usepackage{cite}

\usepackage{ifthen}
\ifthenelse{\equal{\users}{final-layout}}{}{
\usepackage{fancyhdr}
\pagestyle{fancy}
\headheight=28pt\headwidth=17cm
\definecolor{gray}{gray}{0.5}
\rhead{\color{gray}Swelling in viscoelastodynamics Eulerian\\
T.Roub\'\i\v cek \& U.Stefanelli}
\chead{}
\lhead{Version \version, file: \jobname.tex,
\\
compiled:
\number\day.\number\month.\number\year\ at
\the\hour:\ifnum\minute<10 0\fi\the\minute\ h\ \ \ \ \ 
\underline{\Large\color{red}\WHO's working on!}
}
}
 
\newcount\hour \newcount\minute
\hour=\time
\divide \hour by 60
\minute=\time
\loop \ifnum \minute > 59 \advance \minute by -60 \repeat
 
\usepackage[normalem]{ulem}
 
\usepackage[normalem]{ulem}
\usepackage{ifthen}
\usepackage{color} 

\ifthenelse{\equal{\users}{final-layout}}{
	
	\newcommand{\INSERT}[1]{#1}
	\newcommand{\COMMENT}[1]{}
	\newcommand{\COMMENTGT}[1]{}
	\newcommand{\TODO}[1]{}
	\newcommand{\INTERNAL}[1]{}
	\newcommand{\QUESTION}[1]{}
	\newcommand{\DELETE}[1]{}

	\newcommand{\REM}[1]{\marginpar{\bfseries\tiny{\color{blue}}}}
    \newcommand{\MARGINOTE}[1]{}
}
{
	\newcommand{\INSERT}[1]{{\color{blue}\uuline{#1}\color{black}}}
	\newcommand{\COMMENT}[1]{{\color{red}\uuline{#1}\color{black}}}
	\newcommand{\COMMENTGT}[1]{{\hfill\large\color{red}***{#1}***\color{black}\hfill}\\}
	\newcommand{\TODO}[1]{{\color{red}\uuline{#1}\color{black}}}
	\newcommand{\INTERNAL}[1]{\footnote{#1}}
	\newcommand{\QUESTION}[1]{{\color{brown}\uuline{#1}\color{black}}}
	\newcommand{\DELETE}[1]{{\color{red}\sout{#1}\color{black}}}

	\newcommand{\REM}[1]{\marginpar{\bfseries\tiny{\color{blue}#1}}}
\newcommand{\MARGINOTE}[1]{\marginpar{\color{red}\tiny\texttt{#1}}}
}

\newcommand\DT[1]{\mathchoice
                 {{\buildrel{\hspace*{.1em}\text{\LARGE.}}\over{#1}}}
                 {{\buildrel{\hspace*{.1em}\text{\LARGE.}}\over{#1}}}
                 {{\buildrel{\hspace*{.1em}\text{\Large.}}\over{#1}}}
                 {{\buildrel{\hspace*{.1em}\text{\large.}}\over{#1}}}}
\newcommand\pdt[1]{\frac{\partial{#1}}{\partial t}} 
\newcommand{\lineunder}[2]{\LU{\begin{array}[t]{c}\underbrace{#1}\vspace*{.5em}\end{array}}{\mbox{\footnotesize\rm #2}}}
\newcommand{\LU}[2]{\begin{array}[t]{c}#1\vspace*{-1em}\\_{#2}\end{array}}
\newcommand{\linesunder}[3]{\LSU{\begin{array}[t]{c}\underbrace{#1}\vspace*{.5em}\end{array}}{\mbox{\footnotesize\rm #2}}{\mbox{\footnotesize\rm #3}}}
\newcommand{\LSU}[3]{\begin{array}[t]{c}#1\vspace*{-1em}\\_{#2}\vspace*{-.5em}\\_{#3}\end{array}}

\newcommand{\divS}{\mathrm{div}_{\scriptscriptstyle\textrm{\hspace*{-.1em}S}}^{}}
\newcommand{\nablaS}{\nabla_{\scriptscriptstyle\textrm{\hspace*{-.3em}S}}^{}}

\def\Vdots{\!\mbox{\setlength{\unitlength}{1em}
\begin{picture}(0,0)
\put(-.07,0){.}
\put(-.07,.3){.}
\put(-.07,.6){.}
\end{picture}\hspace*{.2em}}}
\usepackage{mathrsfs}   
\usepackage{eucal}      

  \def\bbI{{\bm{I}}}

\def\FG{\boldsymbol}
   
 \def\ee{{\FG e}} \def\ff{{\FG f}} 
  
\def\jj{{\FG j}}   
 \def\nn{{\FG n}}

\def\vv{{\FG v}} \def\ww{{\FG w}} \def\xx{{\FG x}} 
\def\yy{{\FG y}}  
   
\def\DD{{\FG D}} 
\def\FF{{\FG F}} 
\def\II{{\FG I}} 

\def\MM{{\FG M}}   
   
\def\SS{{\FG S}} \def\TT{{\FG T}} 
  \def\XX{{\FG X}}

\newcommand{\R}{\mathbb R}
\newcommand{\N}{\mathbb N}
\newcommand{\Nabla}{{\nabla}}
\newcommand{\Sv}{{\bm D}}
\newcommand{\Se}{{\bm S}}

\newcommand{\Fe}{\FF_{\hspace*{-.2em}\mathrm e^{^{^{}}}}}

\newcommand{\Fes}{\FF}
\newcommand{\Fs}{\FF_{\hspace*{-.2em}\mathrm s^{^{^{}}}}}

\newcommand{\EE}{{\bm e}}

\newcommand{\pl}{\partial}
\newcommand{\eq}[1]{(\ref{#1})}

\renewcommand{\d}{\mathrm d}  
\newcommand{\barOmega}{\hspace*{.2em}{\overline{\hspace*{-.2em}\varOmega}}}

\newcommand{\bulet}{\text{\footnotesize$\,\bullet\,$}}

\newtheorem{theorem}{Theorem}[section]
\newtheorem{lemma}[theorem]{Lemma}
\newtheorem{definition}[theorem]{Definition}

\newtheorem{proposition}[theorem]{Proposition}

\numberwithin{equation}{section}

\newcommand{\UUU}{\color{black}}
\newcommand{\MMM}{\color{black}}
\newcommand{\TTT}{\color{black}}
\newcommand{\EEE}{\color{black}}
\def\OK{\color{black}}

\begin{document}

\allowdisplaybreaks
 
\bigskip\bigskip\bigskip

\noindent{\Large\bf 
Viscoelastodynamics of swelling porous solids\\[.2em]at large strains by an
Eulerian approach
}

\bigskip\bigskip\bigskip\bigskip

\noindent{\bf Tom\'a\v s Roub\'\i\v cek}\footnote{Mathematical Institute, Charles University,
Sokolovsk\'a 83, CZ-186~75~Praha~8,  Czech Republic,\\\hspace*{1.6em} email: ${\texttt{tomas.roubicek@mff.cuni.cz}}$}
\footnote{Institute of Thermomechanics, Czech Academy of Sciences,
Dolej\v skova 5, CZ-18200~Praha~8, Czech Rep.}
\& {\bf Ulisse Stefanelli}\footnote{Faculty of Mathematics, University of
  Vienna, Oskar-Morgenstern-Platz 1, 1090 Vienna, Austria and Vienna Research Platform on Accelerating Photoreaction Discovery, University of Vienna, W\"ahringerstrasse 17, A-1090 Vienna, Austria\\\hspace*{1.6em} email: ${\texttt{ulisse.stefanelli@univie.ac.at}}$}
\footnote{Istituto di Matematica Applicata e Tecnologie Informatiche
  {\it E. Magenes} - CNR, v. Ferrata 1, 27100 Pavia, Italy}


\bigskip\bigskip\bigskip

{\small

  \noindent{\bf Abstract}.
A model
  of saturated hyperelastic
  porous solids at large strains
 is formulated and analysed}. The material response is assumed to be
of a viscoelastic Kelvin-Voigt type and
  inertial effects are considered, too.
  The flow of the diffusant is driven by the
gradient of the chemical potential and is coupled to the mechanics via
the occurrence of swelling and squeezing. Buoyancy effects due to the evolving mass
density in a gravity field are covered. Higher-order viscosity is also included,
allowing for physically relevant stored energies and
local invertibility of the deformation.  
The whole system is formulated in
a fully Eulerian form in terms of rates.
The energetics of the
model is discussed and the existence and regularity of weak solutions is proved by a combined regularization-Galerkin
approximation argument.

\medskip

\medskip


\bigskip\bigskip

\baselineskip=16pt

\def\TRACTION{\bm{f}}
\def\GRAVITY{\bm{g}}
\def\rhoR{\varrho_\text{\sc r}}
\def\RRR{\text{\sc r}}
\def\M{m}
\def\MM{M}

\section{Introduction} 
The  \EEE poromechanics of deformable 
media is a classical part of 
continuum mechanics
of solids, bordering with fluid-solid mechanics and mixtures, and 
having a vast application ground, from petroleum engineering, to geology, soil and rock
mechanics, polymers, etc. Correspondingly, literature is abundant,
see, e.g., the \UUU monographs  \EEE 
 \cite{Boer05TCMP,Cush97PFHP,Stra17MAMP},  
and a whole hierarchy of models
 is available \EEE
\cite{Cush97PFHP,Raja07HAMF}, 
tailored to the description of different aspects at different scales. \EEE

\OK In this paper we focus on a simple \EEE  
phenomenological model \OK describing \EEE  
\OK the flow of a diffusant in a \EEE saturated poroelastic
permeable \OK medium under the assumption that  the diffusant flux is
\EEE governed by \OK Fick's law (here, to be possibly referred as
Darcy's law, for the setting is purely mechanical). \EEE  
\OK We assume the systems to be isothermal and consider the coupled
evolution of the solid and of a single fluid, whose content plays the
role of \EEE 
an internal variable. \OK Indeed, the case of
the \EEE
phenomenological model \OK under scrutiny here  is
\OK justified \EEE  if, among many other simplifications, the flow of a fluid through
the solid is quasistatic and sufficiently slow (i.e., in particular no
\OK inertial effects in the diffusant are considered) \EEE  
and viscous effects \OK in \EEE the fluid  \OK are \EEE neglected.

We focus \OK on the case of possible large strains \EEE and use the
{\it multiplicative decomposition} of the total deformation gradient
\OK into \EEE an elastic strain and a swelling distortion. \OK This
is indeed a classical assumption, especially in connection with 
swelling in \EEE soft materials (like gels) under large strains, see, e.g.,
\cite{BaeSri04DFES,CheAna10CTFP,CheAna11TMCT,CuGaTe17SGCF,DLReAn14CHTP,DroChr13CEFE,DuSoFi10TSMF,LuNaTe13TASI}
or \cite{Anan12CHTT} for the \EEE 
coupling with
inelastic strain (reflecting plasticity or creep).  
The swelling distortion can be modeled, alternatively to the
multiplicative decomposition of the total strain, by
adopting the Biot model \cite{Biot41GTTS}
for the large strains, as used in \cite{HonWan13PFMS,RohLuk17MLDF}
for neo-Hookean material.

The liquid content may influence not only the stress-free configuration through
the mentioned swelling distortion but also the elastic response. \OK
In particular, this amounts in modeling \EEE  
{\it elastic softening} \OK effects. An everyday example of such
phenomenon is the  
soaking of dried \EEE legumes  
\OK which exhibit remarkable swelling accompanied with elastic softening
\OK by increasing \EEE wetting.
 
\OK A specific feature of the model is that it is fully \EEE {\it
  Eulerian}, i.e., \OK it is \EEE
formulated in {\it actual coordinates} instead of a \OK referential
ones. \OK In the frame of the analysis of hyperelastic solid response this is
not common. Still, it allows for some simplification, for it avoids
the need for implementing pullback/pushforward of fileds from the
reference to the actual configuration, ultimately simplifying  transport
coefficients. Let us note however,  that alternative \EEE 
\OK Lagrangian  formulations have been considered in  
\cite{Roub21CHEC,RouSte18TEPR,RouTom18TMPE} or
\cite[Sect.9.6]{KruRou19MMCM}. \OK In spite of the above mentioned
specific analytic intricacies, in contrast with our current Eulerian
one\INSERT{,} these Lagrangian models allow for a possible treatment of
nonhomogeneous Dirichlet conditions on the solid. \EEE 

\OK Additional remarkable features of our model are its fully dynamical nature,
including the description of inertial forces \EEE (thus allowing for elastic
wave propagation), \OK the \UUU attainment \EEE of \EEE local noninterpenetrability
(in the sense that the deformation
gradient is invertible everywhere), and \OK the possibility of
considering \EEE physically relevant stored \OK energies \EEE
(i.e., nonconvex and not necessarily bounded for degenerating
Jacobian of \OK the \EEE deformation).

\OK Extensions to  
multi-porosity or multi-component flows, \OK as well as combinations
\EEE  
with
\OK additional \EEE processes \OK featuring \EEE other evolving internal
variables (as porosity or damage or an inelastic strain) is possible but not
considered here. In addition, one could include thermal effects by \EEE
considering also heat generation and transfer,
possibly with phase transitions. \INSERT{For a metal-hydrid 
phase transition (coupled possibly with magnetic
effects and ferro-to-paramagnetic phase transformation) within
hydrogen diffusion in metals see $\text{\cite{RouTom18TMPE}}$.}

The model is formulated and its energetics is presented in
Section~\ref{sec-model}. The existence of weak solutions is then shown
in Section~\ref{sec-anal}. Here, we follow a 
regularization and  
\OK Galerkin \EEE approximation \OK strategy. This is \EEE combined
with \OK trasport \EEE theory  
by \OK a \EEE regular velocity field.  

\section{The model}\label{sec-model}

 We devote this section to the presentation of the model and its
energetics. After some preparation, these are to be found in
Subsection \ref{sec:model} and \ref{sec:energy-balance}, respectively.

\OK Before going on, let us introduce the \EEE main notation used in this
paper, as in the following table:

\begin{center}
\fbox{
\begin{minipage}[t]{0.39\linewidth}\small\smallskip
$\vv$ velocity (in m/s),\\
$\varrho$ mass density (in kg/m$^3$),\\
$\FF=\Fe\Fs$ deformation gradient,\\
$\Fe$ elastic strain,\\
$\Fs=\lambda(z)\bbI$ swelling distortion,\\
$z$ diffusant content\\
$\bm{T}_{\rm tot}=\TT
{+}\DD$ Cauchy stress (in Pa)\\
$\TT$ elastic (conservative) stress \\
$\DD$ viscous (dissipative) stress 
\smallskip \end{minipage}
\begin{minipage}[t]{0.49\linewidth}\small\smallskip
$\lambda(z)$ swelling stretch\\
$\varphi=\varphi(\Fe,z)$ stored energy (in J/m$^3$=Pa),\\
$\mu$ chemical potential (pore pressure, in Pa)\\
$\M$ mobility (diffusion) coefficient (in 
m$^3$s/kg)\\ 
$\ee(\vv)=\frac12\Nabla\vv^\top\!+\frac12\Nabla\vv$ small strain rate (in s$^{-1}$),\\
$\zeta=\zeta(z,\cdot)$ viscosity dissipation potential,\\
$\GRAVITY$ external load (gravity acceleration in m/s$^{2}$),\\
$\TRACTION$ traction load (in N/m$^2$),\\
$(\cdot)\!\DT{^{}}$ convective derivative
\smallskip \end{minipage}
}\end{center}

\vspace{-.9em}

\begin{center}
{\small\sl Table\,1.\ }
{\small
Summary of the basic notation.
}
\end{center}

\subsection{Geometric preliminaries}
Let us start by recalling some basic notion from the general
theory of large deformations in continuum mechanics. Note that we
limit ourselves in introducing some minimal frame, to serve the sole
purpose of presenting the model. In particular,
no completeness is claimed and we refer the reader. e.g., to  
\cite{GuFrAn10MTC,Mart19PCM} for additional material.   

Assume to be given the {\it deformation} $\yy:I{\times}\varOmega\to\R^d$,
$d=1,2,3$, where $I=[0,T]$ and $T>0$ is some final time. For all given times
$t\in I$, the deformation maps the  {\it reference configuration}
$\varOmega\subset\R^d$ of the deformable body to its {\it actual} configuration
$\yy(t,\varOmega)$, a subset of the {\it physical space} $\R^d$. In what
follows, we indicate {\it referential coordinates} by $\XX\in\varOmega$
and {\it actual coordinates} by $\xx\in\R^d$. By assuming 
$\yy(t,\cdot)$ to be globally invertible, we indicate the inverse by
$\bm{\xi}(t,\cdot)=\yy^{-1}(t,\cdot):\yy(t,\varOmega)\to\varOmega$\TTT; 
standardly $\bm\xi$ is called the {\it return} (or the reference) {\it mapping}
or sometimes {\it inverse motion}. \EEE

Let $Q$ indicate any physical quantity (scalar, vectorial,
tensorial), supposed to be attached to a specific point  $\xx$
of the deformed body at a specific time $t$. The quantity $Q$ can be
expressed in referential coordinates as $Q_\RRR(t,\XX)$ by letting
\begin{equation}
Q_\RRR(t,\XX)= Q(t,\yy(t,\XX)). \label{eq:la}
\end{equation}
Equivalently,
given any physical quantity $Q_\RRR$ (scalar, vectorial,
tensorial), supposed to be attached to a specific referential position
$\XX$ at a specific time $t$ one can express it in actual
coordinates as $Q(t,\xx)$ by posing 
\begin{equation}  Q(t,\xx) = Q_\RRR (t,\bm{\xi}(t,\xx)).
  \label{eq:eu}
\end{equation}
We call $Q$ and $Q_\RRR$ the {\it
  Eulerian} and the {\it Lagrangian or referential} representations of
the quantity, respectively.

Given $\yy=(y_1,\dots y_d)$, we define the
{\it deformation gradient} $\FF_\RRR$ and the {\it referential velocity}
$\vv_\RRR$ as
$$(\FF_\RRR (t,\XX))_{iK}^{}=\frac{\partial y_i}{\partial X_K}(t,\XX) \
\ \ \text{and } \ \ \vv_\RRR(t,\XX) = \frac{\rm d}{{\rm d} t}
\yy(t,\XX)  $$
for indices running from $1$ to $d$. \OK Here and in the following we indicate with ${\rm
  d}/{\rm d} t$ the derivative with respect to time of a time
dependent function, as opposed to the symbol $\pdt{}$ which denotes
the partial time derivative. In the specific case of
$\vv_\RRR(t,\XX)$ these two derivatives obviously coincide. \EEE
The corresponding Eulerian representations from \eqref{eq:eu} are
\begin{align}\label{composition}
\OK \FF(t,\xx) =\FF_\RRR (t, \bm{\xi} (t,\xx)) \ \ \text{and} \ \
\vv(t,\xx) =\vv_\RRR (t, \bm{\xi} (t,\xx))\,. \EEE
\end{align} 
The Eulerian velocity $\vv$ is then used to define the
{\it material derivative} $\DT q(t,\xx)$ of any scalar Eulerian quantity
$q(t,\xx)$ as
$$
\DT q(t,\xx) = \frac{\partial }{\partial t}q(t,\xx)
+ \nabla q(t,\xx){\cdot}\vv(t,\xx) =
\Big( \frac{\partial}{\partial t}+\big(\vv(t,\xx){\cdot}\nabla\big)\Big)
q(t,\xx)\,,
$$
where the differentiation $\nabla$ is, of course, taken with respect to
actual coordinates. Similarly, one defines the material derivative of
a vectorial or tensorial quantity by arguing on coordinates. In
particular, this allows us to check that
\begin{equation}
   \frac{\rm d}{{\rm d} t} Q_\RRR (t,\XX) \stackrel{\eqref{eq:la}}{=}
\frac{\rm d}{{\rm d} t} Q(t,\yy(t,\XX)) = \DT Q(t,\xx)\label{eq:equi}
\end{equation}
For any sufficiently smooth quantity $Q$. In particular, we have that
\begin{equation}
 \label{eq:xi0}
\DT{\bm\xi}(t,\xx) = \frac{\rm d}{{\rm d} t}\bm\xi(t,\yy(t,\XX))
=  \frac{\rm d}{{\rm d} t} \XX={\bm0}.
\end{equation}
Note that property \eqref{eq:equi} in particular implies the product rule 
$ 
\UUU (Q_1 Q_2)^{\text{\LARGE.}}
= \DT Q_1 Q_2 + Q_1 \DT Q_2$.

By applying the classical chain rule we get that
\begin{align*}
  &\frac{\rm d}{{\rm d} t} (\FF_\RRR (t,\XX))_{iK} =  \frac{\rm
  d}{{\rm d} t}\frac{\partial  y_i}{\partial X_K} (t,\XX) =
  \frac{\partial  }{\partial X_K} \frac{\partial  y_i}{\partial t}(t,\XX) \\
  &\quad = \frac{\partial v_{\RRR i} }{\partial X_K}
   (t,\XX)  \stackrel{\eqref{eq:la}}=
   \frac{\partial v_i }{\partial X_K}(t,\yy(t,\XX))  \\
  &\quad =   \frac{\partial  v_i}{\partial x_j}(t,\yy(t,\XX))
  \frac{\partial  y_j}{\partial X_K}(t,\XX)
  =(\nabla \vv(t,\xx))_{ij}^{}\FF_\RRR (t,\XX)_{jK}^{} 
\end{align*}
where, here and below, we use the summation convention over repeated
indices.
Owing to relation \eqref{eq:equi}, the latter reads in Eulerian
coordinates as
\begin{equation}  \DT \FF(t,\xx) = \nabla \vv(t,\xx) \,
\FF(t,\xx).\label{eq:FF}
\end{equation}
\UUU From here on, we formulate the model in terms of the velocity
$\vv$ and the deformation gradient $\FF$ only, without explicit
reference to the deformation $\yy$. \MMM Note that $\yy$ can be
reconstructed by taking the inverse of \UUU 
\MMM $\bm{\xi}$, which is solving equation \eqref{eq:xi0} and is at
least {\it locally} injective. \UUU
The geometric relation
\eqref{eq:FF} will then guarantee that $\FF = \nabla \yy$ for the
reconstructed deformation. \MMM Both \UUU operations are admissible
in our regularity frame, see Definition 3.1 later on, \MMM hinging
on the invertibility of the return mapping $\bm{\xi}$. In
case  $\bm{\xi}$ happens to be {\it globally} injective, the reconstruction of
$\vv$ and $\FF$ can be globally performed. Note however that such
global injectivity is not granted by the model, given the assumed
boundary conditions, which are not fixing tangential deformations. On
the other hand, if boundary deformation were fixed and invertible, we could resort to
the classical theory in \cite{Ball,Kroemer} and deduce global injectivity. \EEE

By using the elementary identity ${\bm0} = \frac{\d}{\d t}
(\FF_\RRR\FF_\RRR^{-1}) =  (\frac{\d}{\d t}\FF_\RRR) \FF_\RRR^{-1} +
\FF_\RRR (\frac{\d}{\d t}\FF_\RRR^{-1})  $ and relation
\eqref{eq:FF}, we also get
\begin{equation}
\UUU(\FF^{-1}(t,\xx))^\text{\LARGE.}\EEE
= - \FF^{-1}(t,\xx)
  \nabla \vv(t,\xx)\,.\label{eq:FF2}
\end{equation}

We now use equations \eqref{eq:FF} and \eqref{eq:FF2} with 
equivalence \eqref{eq:equi} and Jacobi's formula $\frac{\d}{\d t}\det A(t) =\det
A(t) \, \text{tr} (A^{-1}(t)\frac{\d}{\d t}A(t))$, valid for any
sufficiently smooth map $t \mapsto A(t)$ with $A(t)$ invertible, in order to get
\begin{align}\label{DT-det}
   \UUU(\det \FF )^\text{\LARGE.}\EEE
  & =  \frac{\rm
  d}{{\rm d} t} \det \FF_\RRR =  (\det \FF_\RRR )\,\text{tr}\Big(\FF_\RRR^{-1}
  \frac{\d}{\d t}{\FF_\RRR}\Big)\\\nonumber 
  &\!\!\stackrel{\eqref{eq:FF}}{=} (\det \FF  )\, \text{tr} \big(\FF
  ^{-1} \nabla \vv \FF \big) =(\det \FF )\, {\rm div} \vv\,.
\end{align}
Moving from the latter, we also get that
\begin{align}
   &
   \UUU\left( \frac{1}{\det \FF}\right)^\text{\LARGE.}\EEE
   =  \frac{\rm d}{{\rm d} t} \Big( \frac{1}{\det \FF_\RRR}\Big)=  -\frac{(\det \FF_\RRR )\, \text{tr} \left(\FF_\RRR^{-1}\frac{\rm
  d}{{\rm d} t}  \FF_\RRR\right)}{(\det \FF_\RRR)^2} = - \frac{  {\rm
     div} \,\vv}{\det \FF}.\label{DT-det2}
\end{align}

\subsection{The governing equations}

The state of the deformable body undergoing deformation and swelling is
classically described in terms of the actual density $\varrho(t,\xx)$,
the deformation $\yy(t,\xx)$, and the scalar variable $z(t,\xx)$
expressing the pointwise solvent content in Eulerian coordinates. The
evolution of the body is then described by the system
\begin{subequations}\label{eq:model}
  \begin{align}
    &
    \pdt\varrho + \text{div}\, (\varrho \vv)
      =0,\label{eq:model1}\\
    & \varrho \DT{\vv} -\text{div}\, {\bm{T}}_{\rm tot} = \varrho
      \bm{g}, \label{eq:model2}\\[2mm]
    & \DT z -  \text{div}\, (m \nabla \mu)=0. \label{eq:model3}
  \end{align}
\end{subequations}

Here,  relations \eqref{eq:model1}-\eqref{eq:model2} are the
classical conservation of mass and momentum, where ${\bm{T}}_{\rm tot}$
represents the total {\it Cauchy stress} and $\bm{g}$ is the
{\it gravity accelleration}. The kinetic relation
\eqref{eq:model3} describes the transport and diffusion
of the solvent content, in
dependence of the gradient of the {\it chemical potential} $\mu$,
which is additionally  
modulated by the {\it mobility} coefficient $m$. 
Costitutive choices for the quantities $\bm{T}_{\rm tot}$, $\mu$, and
  $m$ are made in Subsection \ref{sec:constitutive} below.

Relations \eqref{eq:model} are to be fulfilled in the deformed domain
$\yy(t,\varOmega)$ for $t\in I$ and have to be complemented by initial and
boundary conditions, see Subsection \ref{sec:model} below. Let us
anticipate that we impose the {\it impenetrability condition} $\vv{\cdot}\bm{n}=0$,
where $\bm{n}$ represents the outward unit normal at the boundary of
the deformed domain. Note that this condition, although possibly being
restrictive with respect to some
applications,  greatly expedits the analysis, for it guarantees that 
$\yy(t,\varOmega)\equiv \varOmega$ for all times $t\in I$. In
particular, one  is actually  asked to solve \eqref{eq:model} on the
cylinder $I{\times}\varOmega$.

Before moving on,
  let us observe that the mass balance \eqref{eq:model1}  
can be equivalently rewritten as $\DT \varrho + \varrho \,
\text{\rm div}\, \vv = 0 $. One can hence use \eqref{DT-det} in order to compute
\begin{align*}
\UUU (\varrho \det \FF)^{\text{\LARGE.}} \EEE
= \DT \varrho \det \FF + \varrho \UUU (\det \FF)^{\text{\LARGE.}}\EEE
\stackrel{\eqref{DT-det}}{=} (\DT \varrho +\varrho \,
\text{\rm div}\, \vv  ) \, \det\FF=0.
\end{align*}
This in particular entails that $\varrho_\RRR (\cdot,\XX) \EEE \det \FF_\RRR (\cdot,\XX) $ is
constant in time for all $\XX\in \varOmega$. Hence,
$$
\varrho_\RRR(t,\OK\XX\EEE)\det\FF_\RRR(t,\OK\XX\EEE)=\varrho_\RRR (0,\XX) \det\FF_\RRR (0,\XX) \,.
$$ 
Passing to Eulerian variables  the latter gives
$$ \varrho (t,\xx) \det \FF(t,\xx) =\varrho (0,\xx) \EEE
\det \OK \FF (0,\xx) \EEE.$$
In particular, provided that relation \eq{eq:FF} holds one can equivalently replace
the continuity equation \eq{eq:model1} and the initial condition
$\varrho \OK (0,\xx) \EEE = \varrho_0 (\xx) $ by 
\begin{align}\label{rho=rho0/detF}
\varrho (t,\xx) =\frac{\varrho_0 (\xx) \det
  \FF (0,\xx) \EEE}{\det\FF (t,\xx) }.
\end{align}
For the sake of later use, we define $\varrho_\RRR(\xx) : =\varrho_0 (\xx) \det
  \FF (0,\xx) \EEE$, which is given in terms of initial data only.

\subsection{Constitutive relations}\label{sec:constitutive}

Let us now fix our \UUU constitutive \EEE choices in relations
\eqref{eq:model1}-\eqref{eq:model3}, leading to the final formulation
of our model in \eqref{Euler-large-diff},
below.

We start by classically assuming that the deformation strain can be multiplicative decomposed as
\begin{align}\label{KLL}
  \FF=
\Fe\Fs.
\end{align}
Here, $\Fe$ denotes the elastic strain whereas $\Fs$ is the strain
associated with swelling. As swelling effects are usually assumed to be
purely volumetric and isotropic, we let
\begin{equation}
  \Fs=\lambda(z)\II
\label{eq:swell}\end{equation}
where the smooth scalar \TTT ``swelling'' \EEE
function $\lambda:[0,1] \to (0,+\infty)$
indicates the stress-free reference volume at solvent-content level $z$ and $\II$ is the identity second-order tensor.

We also assume the total stress $\bm{T}_{\rm tot}$ to be additively decomposed as
\begin{equation}
\bm{T}_{\rm tot} = \bm{T} + \bm{D}\,,\label{eq:TTdec}
\end{equation}
where $\bm{D}$ and $\bm{T}$ denote the viscous (dissipative) and the
inviscid (conservative) stresses, respectively.   

In order to specify constitutive relations, we introduce the stored
energy \TTT in the actual configuration \EEE 
\begin{equation}\label{phi}
  (\Fe,z)\mapsto \int_\varOmega 
\varphi(\Fe,z)
+\delta_{[0,1]}^{}(z) \, \d\xx.
\end{equation}
\MMM By using \eq{KLL} with \eq{eq:swell} so that $\Fe=\FF/\lambda(z)$
one can rewrite equivalently the stored energy in terms of $\FF$ as \TTT
$$(\FF,z)\mapsto \int_\varOmega 
\varphi\Big(\frac\FF{\lambda(z)},z\Big)
+\delta_{[0,1]}^{}(z)\,\d\xx
\,.$$\EEE  
Here, $\varphi$ is the \UUU hyperelastic \EEE energy density \TTT
and $\delta_{[0,1]} :\R\to\{0,+\infty\}$ denotes the indicator
function of the interval $[0,1]$ (namely, $\delta_{[0,1]}(z)=0$
if $z \in [0,1]$ and $\delta_{[0,1]}(z)=+\infty$ otherwise), which in
particular forces 
$z$ to take value in $[0,1]$ only.

\MMM We \EEE define the chemical potential $\mu$ by taking \MMM the 
 variation \EEE  of the stored energy with respect to 
 \TTT $z$\EEE, \MMM namely  
 \begin{equation}
   \label{eq:mu}
   \mu\in \varphi_z'\Big(\frac{\Fes}{\lambda(z)},z\Big)- 
\varphi_{\Fe}'\Big(\frac{\Fes}{\lambda(z)},z\Big){:}\Fes\frac{\lambda'(z)}{\lambda^2(z)} + N_{[0,1]}(z).   
 \end{equation}

 Here, primes denote (partial) differentiation and $N_{[0,1]}$ is the
 subdifferential in the sense of convex analysis of $\delta_{[0,1]}$,
 namely the (multivalued) {\it normal cone} to $[0,1]$ given by  
 $N_{[0,1]}(z) = 0$ if $z \in (0,1)$, $N_{[0,1]}(z) =
 [0,+\infty)$ if $z =1$, $N_{[0,1]}(z) =(-\infty,0]$ if $z=0$, and
 $N_{[0,1]}(z) = \emptyset$ if $z \not \in [0,1]$.
 
 \MMM In order to specify the
 \TTT conservative \EEE stress $\bm{T}$ \MMM we start by computing the
 first Piola-Kirchhoff stress $\bm{P}_\RRR(\bm{X})$ taking the
 variation with respect to $\FF_\RRR$ of the stored energy in
 referential variables, namely,
$$(\FF_\RRR,z_\RRR)\mapsto \int_\varOmega 
\varphi\Big(\frac{\FF_\RRR(\bm{X})}{\lambda(z_\RRR(\bm{X}))},z_\RRR(\bm{X})\Big)\,
\det \FF_\RRR(\bm{X})
+\delta_{[0,1]}^{}(z_\RRR(\bm{X}))\,\d \bm{X},
$$
where we have used $z_\RRR(\bm{X})= z(\bm{y}(\bm{X}))$.
We get
\begin{align*}
  \bm{P}_\RRR(\bm{X}) &= \frac1{\lambda(z_\RRR(\bm{X}))} \varphi_{\Fe}'
\Big(\frac{\FF_\RRR(\bm{X})}{\lambda(z_\RRR(\bm{X}))},z_\RRR(\bm{X})\Big) \det
    \FF_\RRR(\bm{X}) \\
  &\quad + \varphi
\Big(\frac{\FF_\RRR(\bm{X})}{\lambda(z_\RRR(\bm{X}))},z_\RRR(\bm{X})\Big)\,
    {\rm Cof}\, \FF_\RRR(\bm{X})
    \end{align*}
where $\varphi_{\Fe}'$ indicates the derivative of
$\varphi$ in its first variable. Use now the classical position
$\bm{T}_\RRR = (\det \FF_\RRR  )^{-1} \bm{P}_\RRR  \FF_\RRR^\top$ to
conclude that
$$\TT_\RRR=
\frac1{\lambda(z_\RRR)}
\varphi_{\Fe}'\Big(\frac{\Fes_\RRR}{\lambda(z_\RRR)},z_\RRR\Big)\Fes_\RRR^\top\!+\varphi\Big(
\frac{\Fes_\RRR}{\lambda(z_\RRR)},z_\RRR\Big)\bbI.$$
In actual variables, the latter reads \EEE
 \begin{equation}
 \TT=
\frac1{\lambda(z)}
\varphi_{\Fe}'\Big(\frac{\Fes}{\lambda(z)},z\Big)\Fes^\top\!+\varphi\Big(
   \frac{\Fes}{\lambda(z)},z\Big)\bbI.  \label{eq:T}\\
 \end{equation}

 The constitutive equation  for $\bm{D}$ can be obtained
 deduced from a $z$-dependent dissipation potential  
 \begin{equation}
\label{zeta}
\vv\mapsto \int_\varOmega
\zeta(z;\EE(\vv))+\frac{\nu}p|\nabla\EE(\vv)|^p \, \d \xx 
\end{equation}
for some given dissipation density $\zeta$ by taking its variation with
respect to $\EE(\vv)$, namely,
\begin{equation}
  \label{eq:D}
  \bm{D} =  \zeta_\ee'(z;\EE(\vv))-{\rm div}( 
\nu|\nabla\EE(\vv)|^{p-2}\nabla\EE(\vv)).
\end{equation}
The occurrence of the higher-order $\nu$-term corresponds to assuming 
  that the body behaves as a so-called 
  {\it nonsimple material}. This follows the theory 
  by E.~Fried and M.~Gurtin \cite{FriGur06TBBC}, as already
  anticipated in the general nonlinear context of 
  multipolar fluids  by
J.~Ne\v cas at al.\ \cite{NeNoSi89GSIC,NeNoSi91GSCI,NecRuz92GSIV}
or solids \cite{Ruzi92MPTM,Silh92MVMS}, as inspired by 
R.~A.~Toupin \cite{Toup62EMCS} and R.~D.~Mindlin \cite{Mind64MSLE}. 
Such higher-order term in the dissipation ensures
that  
$\nabla\vv$ belongs to $L^1_{\rm w*}(I;L^\infty(\varOmega;\R^{d\times
  d}))$ (weakly$*$ measurable), which guarantees the 
Lipschitz continuity of $\vv(t,\cdot)$, \UUU almost everywhere \EEE in time. This will
turn out crucial in many technical points later on,   in particular in
the estimates
\eqref{test-FF}, \eqref{test-FF+}, and \eqref{test-Delta-r}. Note that,
by dropping such regularity requirement, the treatment of the
transport problem becomes nontrivial due to the possible onset of
singularities, whose occurrence in solids may be debatable 
\cite{AlCrMa19LRCE}.

Eventually, we assume the mobility $m$ to be positive function
depending on $\Fe$ and $z$, namely, $m=m(\Fe,z)$.
 
\def\FeR{\FF_{\hspace*{-.2em}\mathrm e,\text{\sc r}}^{^{^{}}}}
 
\subsection{The model}\label{sec:model}

Following the discussion leading to relation \eqref{rho=rho0/detF}, in
the following we equivalently recast the system \eqref{eq:model}
in terms of the variables $(\vv, \FF, z)$ by dropping the mass
conservation equation \eqref{eq:model1} and requiring the geometric
relation \eqref{eq:FF} instead. Taking also the constitutive
relations \eqref{eq:T}-\eqref{eq:mu} and \eqref{eq:D} into account we
get 
\begin{subequations}\label{Euler-large-CH-plast}
\begin{align}\label{Euler1-large-CH-plast}
&\varrho\DT\vv
={\rm div}(\TT
{+}\Sv)
 +\varrho\GRAVITY\,\ \ \ \text{ with }\ \ \varrho=\frac{\rhoR}{\det\Fes},
     \ \ \ \
     \\&\nonumber\hspace*{4em}
\TT=
\frac1{\lambda(z)}
\varphi_{\Fe}'\Big(\frac{\Fes}{\lambda(z)},z\Big)\Fes^\top\!+\varphi\Big(
\frac{\Fes}{\lambda(z)},z\Big)\bbI\,,
\nonumber\\&\hspace*{4em}
\text{and }\ \ \Sv=\zeta_\ee'(z;\EE(\vv))-{\rm div}(
\nu|\nabla\EE(\vv)|^{p-2}\nabla\EE(\vv))\,,
\nonumber
\\
&\DT\Fes=(\Nabla\vv)\Fes
\,,
\\\label{Euler4-large-CH-plast}
&
\DT z={\rm div}\Big(\M\big(\frac{\FF}{\lambda(z)}
,z\big)\nabla\mu\Big)
 \\&\nonumber\hspace*{4em}
\text{ with }\ \ \mu\in
\varphi_z'\Big(\frac{\Fes}{\lambda(z)},z\Big)-
\varphi_{\Fe}'\Big(\frac{\Fes}{\lambda(z)},z\Big){:}\Fes\frac{\lambda'(z)}{\lambda^2(z)}
+N_{[0,1]}^{}(z).
\end{align}\end{subequations}
We complement the system with the boundary conditions
\begin{subequations}\label{Euler-BC}\begin{align}
&\vv{\cdot}\nn=0\,,\ \ \ 
\big(
(\TT{+}\Sv
)\nn-\divS(
\nu|\nabla\EE(\vv)|^{p-2}\nabla\EE(\vv)
\nn)\big)_\text{\sc t}^{}=\TRACTION\,,\ \ \ 
\\&
\Nabla\EE(\vv){:}(\nn{\otimes}\nn)={\bm0}\,,
\ \ \text{ and }\ \
\M(\FF/\lambda(z)
,z)\nabla\mu{\cdot}\nn+\varkappa\mu=h
\,,
\end{align}\end{subequations}
where the $(d{-}1)$-dimensional surface divergence is defined as
\begin{align}\label{def-divS}
\divS={\rm tr}(\nablaS)\ \ \ \text{ with }\ \
\nablaS \bulet=\nabla \bulet-\frac{\partial \bulet}{\partial\nn}\nn\,,
\end{align}
where ${\rm tr}(\cdot)$ is the trace of a
$(d{-}1){\times}(d{-}1)$-matrix and
$\nablaS  $ denotes the surface gradient. 
Let us again remark the crucial role of the impenetrability boundary
condition $\vv{\cdot}\nn=0$, indeed allowing system
\eqref{Euler-large-CH-plast} to be formulated in the fixed set.
 
We introduce the short-hand notation
\begin{align}
&\widehat\varphi(\FF,z)=\varphi\Big(\frac{\Fes}{\lambda(z)},z\Big)
\ \ \ \ \text{ and }\ \ \ \ 
\widehat\M(\FF,z)=\M\Big(\frac{\Fes}{\lambda(z)},z\Big)
\,.
\label{hat}\end{align}
This allows to rewrite \eqref{Euler-large-CH-plast}
in terms of $\FF$ instead of $\Fe$.
Thus \eqref{Euler-large-CH-plast} can equivalently be written in terms of
$(\vv,\FF,z,\mu)$ 
as
\begin{subequations} \label{Euler-large-diff}
\begin{align}\label{Euler1-large-diff}
     &\varrho\DT\vv
={\rm div}(\TT
{+}\Sv)
 +\varrho\GRAVITY\,\ \ \ \text{ with } \ \varrho=\frac{\rhoR}{\det\Fes},
     \ \ \ \
\TT=\widehat\varphi_{\FF}'(\FF,z)\Fes^\top\!+\widehat\varphi(\FF,z)\bbI\,,
\\&\hspace*{11.2em}
\text{and }\ \ \Sv=\zeta_\ee'(z;\EE(\vv))-{\rm div}(
\nu|\nabla\EE(\vv)|^{p-2}\nabla\EE(\vv))\,,
\nonumber
\\\label{Euler2-diff}
&\DT\Fes=(\Nabla\vv)\Fes
\,,
\\
&
\DT z={\rm div}
\big(
\widehat\M(\FF,z)\nabla\mu\big)
\ \ \ \text{ with }\ \
\mu\in
\widehat\varphi_z'(\FF,z)+
N_{[0,1]}^{}(z)\,.
\label{Euler4-large-diff}
       \end{align}
     \end{subequations}
The last boundary condition in \eq{Euler-BC} can be rewritten
correspondingly as
$$
\widehat\M(\FF,z)\nabla\mu{\cdot}\nn+\varkappa\mu=h.
$$

\subsection{Energy balance}\label{sec:energy-balance}
Let us present the energy balance \UUU underlying \EEE system
\eqref{Euler-large-diff} by testing the three equations respectively by $\vv$,
$\Se
=\widehat\varphi_\FF'(\FF,z)$,  
and $\mu$ and adding up.
After integrating by parts using $\vv{\cdot}\nn=0$, The terms $-{\rm div}\, \TT
{\cdot}\vv$
and $\mu\DT z$ are to be treated jointly   as follows
\begin{align}\label{test-T-nabla-v}
  &  \int_\varOmega
\TT{:}\Nabla\vv+\mu\DT z\,\d \xx=
\int_\varOmega\Se{:}(\Nabla\vv)\Fes+
\widehat\varphi(\FF,z){\rm div}\,\vv
+\mu\DT z\, \d \xx
\\[-.4em]
&\qquad\stackrel{\eqref{Euler2-diff}}{=}\int_\varOmega\!\widehat\varphi_\FF'(\FF,z){:}\DT\FF
+
\widehat\varphi(\FF,z){\rm div}\,\vv
+\widehat\varphi_z'(\FF,z)
\DT z
\, \d \xx
\nonumber\\
&\qquad=\frac{\d}{\d t}\int_\varOmega\widehat\varphi(\FF,z)
\, \d \xx
+\!\!\!\lineunder{\int_\varOmega\nabla\widehat\varphi(\FF,z)
{\cdot}\vv+
\widehat\varphi(\FF,z){\rm div}\,\vv\, \d \xx
}{$=\int_\varGamma
\widehat\varphi(\FF,z)\vv{\cdot}\nn\,d S=0$}\!\!\!\!
\,.
\nonumber\end{align} 

The dissipative terms ${\rm div}\,\Sv{\cdot}\vv$
and $
\mu\,{\rm div}\big(\widehat\M(\FF,z)\nabla\mu\big)$, resulting
by testing \eq{Euler4-large-diff} by $\vv$ and 
\eq{Euler4-large-diff} by $\mu$, can be treated by using twice the Green formula
over $\varOmega$ and once a surface Green formula over $\varGamma$. 
Specifically, using the short-hand notation
$\mathfrak{H}=\nu|\nabla\EE(\vv)|^{p-2}\nabla\EE(\vv)$, we have
\begin{align}\label{Euler-swelling-large-dissip}
&\int_\varOmega{\rm div}\,\Sv{\cdot}\vv+\DT z\mu\,\d \xx=
\int_\varOmega{\rm div}\Big(\zeta_\ee'(z;\EE(\vv))-{\rm div}\mathfrak{H}\Big)
{\cdot}\vv+\mu\,{\rm div}\big(
\widehat\M(\FF,z)\nabla\mu\big)\, \d \xx
\\&\nonumber=\int_\varGamma\Big(\zeta_\ee'(z;\EE(\vv))
-{\rm div}
\mathfrak{H}\Big){:}(\vv{\otimes}\nn)
+
\widehat\M(\FF,z)\nabla\mu{\cdot}\nn\,\d S
\\&\nonumber\qquad-
\int_\varOmega
\Big(\zeta_\ee'(z;\EE(\vv))-{\rm div}\mathfrak{H}\Big){:}\EE(\vv)
+
\widehat\M(\FF,z)|\nabla\mu|^2
\, \d \xx
 \\&\nonumber=
\int_\varGamma\mathfrak{H}{:}(\nn{\otimes}\nn)+
\Big(\zeta_\ee'(z;\EE(\vv))-{\rm div}\,\mathfrak{H}\nn
-\divS(\nn{\cdot}\mathfrak{H})\Big){\cdot}\vv
+
\widehat\M(\FF,z)\nabla\mu{\cdot}\nn
\,\d S
\\&\qquad-
\int_\varOmega\zeta_\ee'(z;\EE(\vv)){:}\EE(\vv)
+\nu|\Nabla\EE(\vv)|^p
+
\widehat\M(\FF,z)|\nabla\mu|^2\, \d \xx\,,
\nonumber\end{align}
where we also used the decomposition of
  $\nabla\vv=(\nn{\cdot}\nabla\vv)\nn+\nablaS\vv$ into
its normal and tangential parts. 

Since under \eqref{Euler-large-diff} we have mass conservation
\eqref{eq:model1} as well, we can   compute 
\begin{align}
  \pdt{}\Big(\frac\varrho2|\vv|^2\Big)=\varrho\vv{\cdot}\pdt\vv
  +\pdt\varrho\frac{|\vv|^2}2
 & =\varrho\vv{\cdot}\pdt\vv
-{\rm div}(\varrho\vv)\frac{|\vv|^2}2
\,.
\label{rate-of-kinetic}\end{align}
By integrating this over $\varOmega$ and using the Green formula and
$\vv{\cdot}\nn=0$, we obtain
\begin{align}
\frac{\d}{\d t}\int_\varOmega\frac\varrho2|\vv|^2\,\d x
=\int_\varOmega
\varrho\vv{\cdot}\pdt\vv
+\varrho\vv{\cdot}(\vv{\cdot}\nabla)\vv\,\d x
-\!\int_\varGamma\varrho|\vv|^2\vv{\cdot}\nn\,\d S
=\int_\varOmega\varrho\DT\vv{\cdot}\vv\,\d x\,.
\label{calculus-convective-in-F}
\end{align}
In particular, the inertial force $\varrho\DT\vv$ tested by $\vv$ can be treated
by \eq{calculus-convective-in-F} while 
the gravity force density $\varrho\GRAVITY$ yields directly the power
of the gravitational field $\varrho\GRAVITY{\cdot} \vv$.

Eventually, one obtains the energy balance
\begin{align}\nonumber
  &\frac{\d}{\d t}
  \int_\varOmega\!\!\linesunder{\frac\varrho2|\vv|^2}{kinetic}{energy}\!\!+\!\!\linesunder{
  \widehat\varphi(\FF,z)}{stored}{energy}\!\!\d \xx
+\int_\varOmega
\!\!\!\!\linesunder{
\zeta_\ee'(z;\EE(\vv)){:}\EE(\vv)
+\nu|\nabla\EE(\vv)|^p
_{_{_{}}}\!}{dissipation rate}{due to viscosity
}\!\!\!\d \xx\\&\hspace{0em}
+\int_\varOmega\!\!\!\linesunder{
\widehat\M(\FF,z)|\nabla\mu|^2}{dissipation rate}{due to diffusion}\!\!\!\d \xx
+\int_\varGamma\!\!\!\!\!\!\!\!\linesunder{\varkappa\mu^2_{_{_{_{_{}}}}}}{dissipation}{rate of influx}\!\!\!\!\!\!\!\!\!\!\d S
=\int_\varOmega\!\!\!\!\!\!\linesunder{\varrho\GRAVITY{\cdot}\vv_{_{_{_{_{}}}}}\!\!\!}{power of}{gravity field}\!\!\!\!\!\!\d \xx+\int_\varGamma\!\!\!\!\linesunder{\TRACTION{\cdot}\vv_{_{_{_{_{}}}}}}{power of\ \ }{traction\ \ }\!\!\!\!\!+\!\!\!\!\!\!\!
\linesunder{h\mu_{_{_{_{_{}}}}}}{chemical}{influx}
\!\!\d S\,.
  \label{Euler-swelling-large-energy}\end{align}

\section{The analysis by Faedo-Galerkin semi-discretization}
\label{sec-anal}

We consider  the Cauchy   problem for the  system
\eq{Euler-large-CH-plast} with boundary conditions  \eq{Euler-BC}. For this, we prescribe the
initial conditions
\begin{align}\label{IC}
\vv|_{t=0}^{}=\vv_0\,,\ \ \ \ \FF|_{t=0}^{}=\FF_0\,,\ \ \text{ and }\ \
z|_{t=0}^{}=z_0\,.
\end{align}

\UUU In the following, we assume $\varOmega\subset \R^n$ to be a nonempty, open,
bounded, connected set with Lipschitz boundary $\varGamma:=\partial
\varOmega$. \EEE
We will use the  following  standard notation  for  Lebesgue and Sobolev
spaces. Namely,
$L^p(\varOmega;\R^n)$ denotes the Banach space of Lebesgue measurable functions
$\varOmega\to\R^n$ whose  $p$-power  of the  Euclidean norm is
integrable and
$W^{k,p}(\varOmega;\R^n)$  is  the space of
$L^p(\varOmega;\R^n)$ functions whose
derivatives of order $k$ are in $L^p(\varOmega;\R^{n\times kd})$.
 We indicate  $W_0^{2,p}(\varOmega;\R^d):=\{\vv\in W^{2,p}(\varOmega;\R^d);\
\vv{\cdot}\nn=0\text{ on }\varGamma\}$  and use the short-hand
notation 
  $H^k=W^{k,2}$.  Given   a Banach space
$X$ and $I=[0,T]$, we use the notation $L^p(I;X)$ for the Bochner
space of Bochner measurable functions $I\to X$ whose norm is in $L^p(I)$, 
and $H^1(I;X)$ for functions $I\to X$ whose distributional derivative
is in $L^2(I;X)$.  The spaces $C(I;X)$ and  $C_{\rm w}(I;X)$
indicate continuous and   weakly continuous
functions $I\to X$, respectively. Dual spaces are denoted by 
$(\cdot)^*$     and
$p'=p/(p{-}1)$  indicates  the conjugate exponent, with the convention
$p'=\infty$ for $p=1$ and $p'=\infty$ for $p=1$.
 For $p<d$,  we indicate  by $p^*$ the exponent from the embedding
$W^{1,p}(\varOmega)\subset L^{p^*}(\varOmega)$, i.e.\
$p^*=pd/(d{-}p)$.  
Occasionally, we will use $L_{\rm w*}^p(I;X)$ for weakly* measurable
functions $I\to X$ for nonseparable spaces $X$ which are duals to some other
Banach spaces (specifically for $L^\infty(\varOmega)$).  

\def\OMEGA{\omega}

The energy balance \eq{Euler-swelling-large-energy} 
delivers formal   a-priori estimates.  Aiming at making this
rigorous, we start by specifying  our assumptions
on the data.  By indicating by ${\rm GL}^+(d)$ the space of $d{\times}d$
matrices with positive determinant, we ask for the following. 
\begin{subequations}\label{ass}
\begin{align}
&\label{ass-phi}
\varphi : \R^{d \times d} \to (-\infty,+\infty], \ \ \varphi  \in C^1(
  {\rm GL}^+(d)  \times\R)\ \ \UUU \exists \MMM \kappa \UUU >0 \ \
  \text{such that} \EEE
\\[2mm]&\nonumber
\ \ \ \ \ \varphi(\Fe,z)\ge \MMM \kappa \UUU/\det\Fe\ \ \text{for \UUU all
         $\Fe$ with \EEE  $\det\Fe>0$}, \\&\nonumber
\ \ \ \ \ \text{$\varphi(\Fe,z)=+\infty$ for $\det\Fe\le0$, and}
\\&\nonumber
\ \ \ \ \ z\mapsto\widehat\varphi(\FF,z)=\varphi\Big(\frac{\FF}{\lambda(z)},z\Big)
\ \ \text{ strongly convex,  uniformly  \UUU w.r.t \EEE $\FF$, \UUU
    namely, \EEE}
\\
&\nonumber
\ \ \ \ \ \UUU 
\forall \FF \in   {\rm GL}^+(d), \ z_0,\,z_1\in\TTT[0,1],\UUU 
\ \theta\in [0,1]\EEE:
\\&\nonumber
\ \ \ \ \ \ \ \ \ \ \widehat\varphi(\FF,\theta z_1+
  (1{-}\theta)z_0) \leq \theta \widehat\varphi(\FF, z_1)+ (1{-}\theta )
  \widehat\varphi(\FF, z_0) - \frac{\kappa}{2}
  |z_1{-}z_0|^2\,,
\\
 &\UUU p>d, \EEE \label{ass-p}
\\&
\zeta:\R\times\R_{\rm sym}^{d\times d}\to\R\ \text{ continuously differentiable},\ \zeta(z,\cdot)\ \text{ convex},
\\&\label{ass-D}
\ \ \ \ \ \exists\MMM\bar \varepsilon \EEE>0\ \forall(z,\ee)\in\R\times\R_{\rm sym}^{d\times d}:\ \ 
\MMM\bar \varepsilon \EEE|\ee|^2\le\zeta(z,\ee)\le(1{+}|\ee|^2)/\MMM\bar \varepsilon \EEE\,,
\\&\lambda\in C^1(\R) \cap W^{1,\infty}(\R)  \ \text{ and }\inf\lambda>0\,,
\\&\M:\R^{d\times d}\times\R\to\R\ \text{ continuous  and  bounded with
 }\ 
\mbox{$\inf_{\R^{d\times d}\times\R}^{}m>0$}\,,
\label{ass-j}
\\&
\label{ass-g}
\GRAVITY\in L^1(I;L^\infty(\varOmega;\R^d))\,,\ \ \ \TRACTION\in L^{p'}(I;L^1(\varGamma;\R^d))\,,
\ \ \ h\in L^2(I;L^{4/3}(\varGamma))\,,
\\&\label{ass-IC}
\vv_0\in L^2(\varOmega;\R^d)\,,\ \ \FF_0\in W^{1,r}(\varOmega;\R^{d\times d})\,,
\ \ z_0\in W^{1,r}(\varOmega)\,,\ \ \rhoR\in W^{1,r}(\varOmega)\,,\ \ r>d\,,
\\&\ \ \ \ \ \ \ \ \ \ \ \ \ \ \ \text{ with} \ \ 
{\rm min}_{\barOmega}^{}\det\FF_0>0\ \ \text{ and }\ \
{\rm min}_{\barOmega}^{}\rhoR^{}>0\,.
\nonumber\end{align}\end{subequations} 
 Note that the concrete form of  $\varphi(\FF,z)$ will 
actually  be relevant only for $z\in[0,1]$.  Still, as in the
proof of
Proposition~\ref{prop-Euler-plast} such constraint is penalized, we
are asked to  define $\varphi$ also for
outside the interval $[0,1]$ in \eqref{ass-phi}, and similarly also
for $\zeta$ and $\M$.

\MMM Before moving on, let us show that conditions \eq{ass-phi} can be
realized in some physically relevant situation. To this aim, let
$\lambda$ be positive.
For
all $\Fe$ with $\det \Fe>0$ let the {\it Ogden-type} energy density be
defined as
$$\varphi(\Fe,z)=f_1(z)g_1(\Fe\Fe^\top) + f_2(z)g_2({\rm Cof}
(\Fe\Fe^\top)) +f_3(z)g_3(\det \Fe)+\frac\kappa{\det\Fe} +h(z)\,,$$
where $h$ is uniformly convex, $g_i,\, f_i $ are smooth, nonnegative, convex,
with $g_i(0)=0$, for $i=1,2,3$. Indeed, the uniform convexity of $z\mapsto \widehat\varphi(\FF,z)$ with
respect to $\FF$ from \eqref{ass-phi} follows from the uniform
convexity of $h$ and from the convexity in $z$ of all other
terms. Such convexity
can be checked by noticing that all such terms have the form
$\eta(\Fe,z)=f(z)g(H(\Fe))   $ where $H$ is a $s$-homogeneous
function. In particular, $$\eta(\Fe,z) = \eta \left(\frac{\FF}{\lambda(z)},z\right) =
f(z)g\left(\frac{H(\FF)}{(\lambda(z))^s} \right).$$
Under general assumptions on $f$ and $g$, convexity can be directly checked by
computing the second derivative with respect to $z$. To simplify
notation, assume $H$ to be scalar valued (which is the case for $i=3$)
and compute
\begin{align*}
  &\eta''_z(\Fe,z) = \eta''_z\left(\frac{\FF}{\lambda(z)},z\right) = f''(z)
  g\left(\frac{H(\FF)}{(\lambda(z))^s} \right) -2s f'(z)  
  g'\left(\frac{H(\FF)}{(\lambda(z))^s} \right)
  \frac{H(\FF)\lambda'(z)}{(\lambda(z))^{s+1}}\\
  &\quad-sf(z) g'\left(\frac{H(\FF)}{(\lambda(z))^s} \right)
  \frac{H(\FF)\lambda''(z)}{(\lambda(z))^{s+1}}+s(s+1)f(z) g'\left(\frac{H(\FF)}{(\lambda(z))^s} \right)
    \frac{H(\FF)(\lambda'(z))^2}{(\lambda(z))^{s+2}} \\
  &\quad+ s^2 f(z) g''\left(\frac{H(\FF)}{(\lambda(z))^s} \right)
    \frac{(H(\FF))^2(\lambda'(z))^2}{(\lambda(z))^{2s+2}}\\
  &\geq s(\lambda(z))^2 g'\left(\frac{H(\FF)}{(\lambda(z))^s} \right)
    \left(\frac{H(\FF)}{(\lambda(z))^s} \right)\left((s{+}1) f(z)  (\lambda'(z)) ^2  {-}2
    f'(z)\lambda'(z) \lambda(z) {-} f(z)  \lambda''(z) \lambda(z)   \right).
\end{align*}
As $g'(r)r \geq g(r) \geq 0$, the latter is nonnegative under
appropriate assumptions on $\lambda$ and $f$, for instance, if  $f$ is
nonincreasing and $\lambda$ is nondecreasing and concave.
 
\EEE

 In order to obtain   a-priori estimates from
\eq{Euler-swelling-large-energy},  a number of technical
points have to be faced. One first 
issue is estimation of the gravity force $\varrho\GRAVITY$  when
tested by the velocity $\vv$, which  can be estimated by
the H\"older/Young inequality as
\begin{align}\label{Euler-est-of-rhs}
\int_\varOmega\varrho\GRAVITY{\cdot}\vv\,\d \xx
&=\int_\varOmega\sqrt{\frac{\rhoR}{\det\FF}}\sqrt{\varrho}\vv{\cdot}\GRAVITY\,\d \xx
\le\Big\|\sqrt{\frac{\rhoR}{\det\FF}}\Big\|_{L^{2}(\varOmega)}
\big\|\sqrt{\varrho}\vv\big\|_{L^2(\varOmega;\R^d)}^{}\big\|\GRAVITY\big\|_{L^\infty(\varOmega;\R^d)}^{}
\\&\nonumber\le
\frac12\bigg(\Big\|\sqrt{\frac{\rhoR}{\det\FF}}\Big\|_{L^2(\varOmega)}^2
\!+\big\|\sqrt{\varrho}\vv\big\|_{L^2(\varOmega;\R^d)}^2\bigg)
\,\big\|\GRAVITY\big\|_{L^\infty(\varOmega;\R^d)}^{}
\\& =  \big\|\GRAVITY\big\|_{L^\infty(\varOmega;\R^d)}^{}
\int_\varOmega
\frac{\rhoR}{2\det\FF\!}
+\frac\varrho2|\vv|^2\,\d \xx\,.
\nonumber\end{align}
The integral on the right-hand side can be then treated by the Gronwall
 lemma, by  relying on the kinetic-energy term and  the
fact that the stored energy controls $1/\det \FF$, i.e., 
 \eq{ass-phi}.  In order to apply the Gronwall lemma one needs
 the qualification \eq{ass-g}  for  
 $\GRAVITY$.
 
 A second  technical issue is estimation of the boundary term
$\ff{\cdot}\vv$,  which will follow along the lines of relation
\eq{est-fv} below, which in turn hinges on a bound on 
$1/\varrho$, cf. \eqref{est+rho}.

Eventually, $\int_\varGamma h\mu\,\d S$ can be estimated by
$\|h\|_{L^{4/3}(\varGamma)}^{}\|\mu\|_{L^4(\varGamma)}^{}\le\|h\|_{L^{4/3}(\varGamma)}^2/\delta+
\delta\|\mu\|_{L^4(\varGamma)}^2\le\|h\|_{L^{4/3}(\varGamma)}^2/\delta+
\delta N\|\nabla\mu\|_{L^2(\varOmega;\R^d)}^2+
\delta N\|\mu\|_{L^2(\varGamma)}^2$, where $N$  indicates the
square of the   norm of
the trace operator $H^1(\varOmega)\to L^4(\varGamma)$.

 Under assumptions \eqref{ass}, the energy balance
\eqref{Euler-swelling-large-energy} thus implies the a-priori estimates
\begin{subequations}\label{est}
\begin{align}\label{est-rv2}
&\|\sqrt{\varrho}\vv\|_{L^\infty(I;L^2(\varOmega;\R^d))}^{}\le C\,,
\\&\label{est-phi}\|\varphi(\Fe,z)\|_{L^\infty(I;L^1(\varOmega))}^{}\le C\,,
\\&\label{est-e(v)}
\|\EE(\vv)\|_{L^2(I;W^{1,p};\R^{d\times d}))}^{}\le C\,,
\\&\label{est-z}
\|z\|_{L^\infty(I{\times}\varOmega
)}^{}\le C\,,\ \ \ \text{ and }
\ \ \
\\&\label{est-mu}\|\mu\|_{L^2(I;H^1(\varOmega))}^{}\le C\,,
\end{align}\end{subequations}
 where, here and in the following, for the sake  of notational
simplicity the symbol $C$ is used to
indicate any positive constant just depending on data and possibly
varying form line to line. In case of need, we will indicate the
dependence of such constant on specific parameters by using
\OK subscripts. \EEE 

 As  $p>d$, estimate  \eq{est-e(v)} prevents the onset  of
singularities  for  the quantities transported  by the
  velocity field.  In particular,  due to qualification of
$\FF_0$ and $\varrho_0=\rhoR/\det\,\FF_0$ in \eq{ass-IC},  
it yields the estimates
\begin{subequations}\label{est+}
\begin{align}
\label{est+Fes}&\|\Fes\|_{L^\infty(I;W^{1,r}(\varOmega;\R^{d\times d}))}\le C_r\,,
\ \ \ \Big\|\frac1{\det\Fes}\Big\|_{L^\infty(I;W^{1,r}(\varOmega))}\le C_r\,,
\\&\label{est+rho}\|\varrho\|_{L^\infty(I;W^{1,r}(\varOmega))}^{}\le C_r\,,
\ \ \text{ and }\ \ \Big\|\frac1\varrho\Big\|_{L^\infty(I;W^{1,r}(\varOmega))}\!\le C_r
\ \ \ \text{ for any $1\le r<+\infty$},
\intertext{ see   the arguments in the proof of
Lemmas~\ref{lem-transport-F+} and \ref{lem2} below. 
From \eq{est-rv2} and \eq{est+rho}, we then  also  have}
&\|\vv\|_{L^\infty(I;L^2(\varOmega;\R^d))}^{}\le
\|\sqrt\varrho\vv\|_{L^\infty(I;L^2(\varOmega;\R^d))}^{}\Big\|\frac1{\sqrt\varrho}\Big\|_{L^\infty(I\times\varOmega)}^{}\le C\,.
                                                                   \label{basic-est-of-v}\end{align}\end{subequations}

 Based on the formal a-priori estimates \eqref{est}-\eqref{est+}
we now specify a notion of weak solution.  
In
particular, we replace the inertial force $\varrho\DT\vv$ in
\eq{Euler1-large-CH-plast} by using  the equality
\begin{align}\label{inertial}
\varrho\DT\vv=\frac{\partial}{\partial t} (\varrho \vv) +
\text{\rm div}\, (\varrho \vv \otimes \vv),
\end{align}
as well as 
$\varrho(0)=\rhoR/\det\FF_0$. Noteworthy, this formula has exploited
the continuity equation \eq{eq:model1}.

\begin{definition}[Weak solutions to \eq{Euler-large-CH-plast}]\label{def}
We call quintuple $(\varrho,\vv,\Fes,
z,\mu)\in L^\infty( I{\times}\varOmega )$ $\times
(L^p(I;W^{2,p}(\varOmega;\R^d))\cap C_{\rm  w }(I;L^2(\varOmega)))\times
L^\infty( I{\times}\varOmega;\R^{d\times d}))\times
C_{\rm w}(I;
L^2(\varOmega))\times L^2(I;H^1(\varOmega))$ a \emph{weak solution}
to the system \emph{\eq{Euler-large-CH-plast}} with initial and boundary
conditions \emph{\eq{Euler-BC}} and \emph{\eq{IC}} if $\vv{\cdot}\nn=0$,
$\det\FF>0$ and $\varrho=\rhoR/\det\FF$ a.e.\ on $I{\times}\varOmega$,
$\varphi(\FF,z),\varphi_z'(\FF,z)\in L^1( I{\times}\varOmega )$,
$\varphi_\FF'(\FF,z)\in L^1( I{\times}\varOmega ;\R^{d\times d})$, $z$ is valued in $[0,1]$, $\vv(0)=\vv_0$, 
\begin{subequations}\label{Euler-weak}\begin{align}
&\label{def-Euler-weak1}
\int_0^T\!\!\!\!\int_\varOmega\bigg(
\Big(
\widehat\varphi_\FF'(\FF,z)\FF^\top
+\zeta_\ee'(z;\EE(\vv))-\varrho\vv{\otimes}\vv\Big){:}\Nabla\widetilde\vv
-\varrho\vv{\cdot}\pdt{\widetilde\vv}
\TTT+\EEE
\widehat\varphi(\FF,z)({\rm div}\,\widetilde\vv)
\\[-.4em]&\hspace{3em}
+\nu|\nabla\EE(\vv)|^{p-2}\nabla\EE(\vv)\Vdots\Nabla\EE(\widetilde\vv)\bigg)
\,\d \xx\d t
=\!
\int_0^T\!\!\!\!\int_\varOmega\!\varrho\GRAVITY{\cdot}\widetilde\vv\,\d \xx\d t
+\int_0^T\!\!\!\!\int_\varGamma\TRACTION{\cdot}\widetilde\vv\,\d S\d t
\nonumber
\intertext{holds for any
$\widetilde\vv\UUU\in C^\infty(I{\times}\MMM \barOmega \EEE;\R^d)$ \EEE with
$\widetilde\vv{\cdot}\nn=0$ and $\widetilde\vv(T)=0=\widetilde\vv(0)$, }
&\int_0^T\!\!\!\!\int_\varOmega\!\Fes{:}\pdt{\widetilde\SS}
+\Big(({\rm div}\,\vv)\Fes{+}(\Nabla\vv)\Fes
\Big){:}\widetilde\SS
+\Fes{:}((\vv{\cdot}\Nabla)\widetilde\SS)\,\d \xx\d t
=-\!\!\int_\varOmega\!
\FF_0{:}\widetilde\SS(0)\,\d \xx\!\label{def-Euler-weak2}
\intertext{holds for any $\widetilde\SS\UUU \in C^\infty(I{\times}
            \MMM \barOmega \EEE;\R^{d\times d})$ \EEE with  $\widetilde\SS(T)={\bm0}$,}
\label{def-Euler-weak3}
&\int_0^T\!\!\!\int_\varOmega
\widehat\M(\FF,z)\nabla\mu{\cdot}\nabla\widetilde z-z\pdt{\widetilde z}
-z{\rm div}(\vv\widetilde z)
\,\d \xx\d t+\int_0^T\!\!\!\int_\varGamma \varkappa\mu\widetilde z\,\d S\d t
\\[-.4em]&\hspace{17em}
=\int_\varOmega\mu_0\widetilde z(0)\,\d \xx
+\int_0^T\!\!\!\int_\varGamma h\widetilde z\,\d S\d t
\nonumber
\intertext{holds for any $\widetilde z \UUU \in C^\infty(I{\times}
            \barOmega)$ \EEE with
$\widetilde z(T)=0$, and
}
&\int_0^T\!\!\!\int_\varOmega\!
\big(\widehat\varphi_z'(\FF,z)
-\mu\big)(\widetilde z{-}z)\,\d \xx\d t\ge0
\label{def-Euler-weak4}\end{align}
\end{subequations}
holds for any $\widetilde z\in L^\infty( I{\times}\varOmega )$ 
valued in $[0,1]$.
\end{definition}

 If the velocity field $\vv$ is in $L^1(I;W^{1,\infty}(\varOmega;\R^d))$,
one classically obtain that regularity of the initial datum is preserved along
the flow \eqref{Euler2-diff}. We provide a rigorous statement in the following
lemma, as well as a proof based on Galerkin approximations. 

\begin{lemma}[Flow of  $\FF$
]\label{lem-transport-F+}
Let $p>d$ and $r>2$. 
Then, for any $\vv\in L^1(I;W^{2,p}(\varOmega;\R^d))$
with 
$\vv{\cdot}\nn=0$  
and any $\Fes_0\in W^{1,r}(\varOmega;\R^{d\times d})$, there exists
 a unique  weak solution $\FF\in C_{\rm w}(I;W^{1,r}(\varOmega;\R^{d\times d}))\cap W^{1,1}(I;L^r(\varOmega;\R^{d\times d}))$
to \eqref{Euler2-diff}   in the sense of 
\eqref{def-Euler-weak2}.  The  estimate
\begin{align}\label{F-evol-est}
\|\FF\|_{L^\infty(I;W^{1,r}(\varOmega;\R^{d\times d}))\,\cap\,
W^{1,1}(I;L^r(\varOmega;\R^{d\times d}))}^{}\le
\mathfrak{C}\Big(\|\Nabla\vv\|_{L^1(I;W^{1,p}(\varOmega;\R^{d\times d}))}^{}\,,\,
\|\Fes_0\|_{W^{1,r}(\varOmega;\R^{d\times d})}^{}\Big)
\end{align}
holds with some $\mathfrak{C}\in C(\R^2)$,   equation 
\eqref{Euler2-diff} holds a.e.\ on $I{\times}\varOmega$,
and $\FF\in C(I{\times}\barOmega;\R^{d\times d})$. Moreover, the mapping
\begin{align}\label{v-mapsto-F}
\vv\mapsto\FF:L^1(I;W^{2,p}(\varOmega;\R^d))\to
L^\infty(I;W^{1,r}(\varOmega;\R^{d\times d}))
\end{align}
is (weak,weak*)-continuous. 
If in addition $\det\Fes_0>0$ on $\barOmega$, then $\det\FF >0 $ on
$ I{\times}\barOmega$;  i.e., $\min_{ I{\times}\barOmega }\det\UUU \FF
\EEE >0$
 uniformly  with respect to bounded velocity fields $\vv$, namely, 
for any  $R>0$   there  exists  $\delta>0$ such that
\begin{align}
\|\Nabla\vv\|_{L^1(I;W^{1,p}(\varOmega;\R^{d\times d}))}^{}\le R\ \ \Rightarrow\ \
\min_{ I{\times}\barOmega }\det\FF\ge\delta\,.
\label{detF>0}\end{align}
\end{lemma}

\begin{proof}
 Let us start by assuming 
$\vv\in L^2(I;W^{2,p}(\varOmega;\R^d))$; the weaker integrability setting of
$\vv\in L^1(I;W^{2,p}(\varOmega;\R^d))$ will be recovered later in the
proof.

Consider the following parabolic regularization of \eqref{Euler2-diff}
\begin{align}\label{F-evol-reg}
\DT\Fes=(\Nabla\vv)\Fes
+\varepsilon{\rm div}(|\nabla\Fes|^{r-2}\nabla\Fes)\,,
\end{align}
 by complementing it by the  additional boundary condition $(\Nabla\Fes)\nn=\bm0$.
 We tackle the  regularized problem \eq{F-evol-reg} by means of a
Faedo-Galerkin approximation. Assume to be given a sequence of nested
finite-dimensional subspaces $\{\TTT U\EEE_k\}_{k\in\N}$ whose union is dense in
$W^{1,\TTT r\EEE}(\varOmega;\R^{d\times d})$. Without loss of generality, we
can  ask for $\Fes_0\in \TTT U\EEE_1$. The classical existence theory for
systems of ordinary differential equations ensures that one can find a solution 
$t \in I \mapsto \FF_k(t) \in \TTT U\EEE_k$ of the Galerkin-approximated
problem for any $k$; more precisely, local in time existence needs to be
combined with maximal prolongation on the whole interval $I$, on the basis
of the  $L^\infty$-estimates below).

Testing  (the Galerkin approximation of)  \eq{F-evol-reg} 
by $\Fes_k$ we can estimate 
\begin{align}\label{test-FF}
&\frac{\d}{\d t}\int_\varOmega\frac12|\Fes_k|^2\,\d \xx
+
\varepsilon\int_\varOmega|\Nabla\Fes_k|^r\,\d \xx
=\int_\varOmega\Big((\nabla\vv)\Fes_k-(\vv{\cdot}\nabla)\Fes_k
\Big){:}\Fes_k
\,\d \xx
\\[-.4em]&
\qquad\ \ 
=\int_\varOmega\!(\nabla\vv)\Fes_k{:}\Fes_k
+\frac{{\rm div}\,\vv}2|\Fes_k|^2
\,\d \xx
\le
\frac32\|\nabla\vv\|_{L^\infty(\varOmega;\R^{d\times d})}^{}
\|\Fes_k\|_{L^2(\varOmega;\R^{d\times d})}^2\,;
\nonumber\end{align}
 where we used the  calculus 
\begin{align}\nonumber
  \int_\varOmega(\vv{\cdot}\nabla)\Fes_k{:}\Fes_k\,\d \xx
  &=\!\int_\varGamma|\Fes_k|^2(\vv{\cdot}\nn)\,\d S
 \\[-.4em]\nonumber&\quad
  -\!\int_\varOmega\!\Fes_k{:}(\vv{\cdot}\nabla)\Fes_k+({\rm div}\,\vv)|\Fes_k|^2\,\d \xx
=-\frac12\int_\varOmega({\rm div}\,\vv)|\Fes_k|^2\,\d \xx,
\end{align}
together with the boundary condition $\vv{\cdot}\nn=0$.  Note in
particular that, in order to perform the latter integration by parts,
the integrability of $\vv$ is required, besides the regularity of 
  $\Nabla\vv$.  
 By the Gronwall inequality 
we obtain the estimate
\begin{align}
\label{Euler-quasistatic-est1-2}
\|\Fes_k\|_{L^\infty(I;L^2(\varOmega;\R^{d\times d}))}\le C\ \ \text{ with }\ \ \|\Nabla\Fes_k\|_{L^r( I{\times}\varOmega ;\R^{d\times d\times d})}\le C\varepsilon^{-1/r}\,.
\end{align}

 At  the Galerkin-discretization level, another legitimate test
  for  \eq{F-evol-reg} is   $\pdt{}\Fes_k$. This allows us to estimate 
\begin{align}\label{test-FF+}
&\int_\varOmega\bigg|\pdt{\Fes_k}\bigg|^2\,\d \xx
+\frac
\varepsilon r\frac{\d}{\d t}\int_\varOmega|\Nabla\Fes_k|^r\,\d \xx
=\int_\varOmega\Big((\nabla\vv)\Fes_k-(\vv{\cdot}\nabla)\Fes_k
\Big){:}\pdt{\Fes_k}\,\d \xx
\\[-.0em]&\nonumber\hspace{6em}
\le
\|\nabla\vv\|_{L^\infty(\varOmega;\R^{d\times d})}^2
\|\Fes_k\|_{L^2(\varOmega;\R^{d\times d})}^2
\\[-.2em]&\hspace{7em}
+C_r\|\vv\|_{L^\infty(\varOmega;\R^d)}^2
\Big(
1+\|\Nabla\Fes_k\|_{L^r(\varOmega;\R^{d\times d})}^r\Big)
+\frac12\bigg\|\pdt{\Fes_k}\bigg\|_{L^2(\varOmega;\R^{d\times d})}^2.
\nonumber\end{align}
Note that we used here that $r>2$. 
As $\vv\in L^2(I;W^{2,p}(\varOmega;\R^d))
\subset  L^2(I;L^\infty(\varOmega;\R^d))$, 
the already obtained estimate  
\eq{Euler-quasistatic-est1-2}, and the Gronwall inequality  imply
that  
\begin{align}
\label{Euler-quasistatic-est2}
&
                                 \Big\|\pdt{\Fes_k}\Big\|_{L^2(
  \varOmega {\times} I ;\R^{d\times d})}\le
                                   C{\rm
                                 e}^{1/(r\varepsilon)} 
\ \ \text{ and }\ \ 
\|\Nabla\Fes_k\|_{L^\infty(I;L^r(\varOmega;\R^{d\times d\times d}))}\le C
{\rm e}^{1/(r\varepsilon)}
\,.
\end{align}

Keeping  $\varepsilon>0$ fixed,  the above  estimates allow \TTT us \EEE to
pass to the limit as $k\to\infty$ by standard arguments for quasilinear
parabolic equations\TTT; realize that all lower-order terms are linear
while the last, highest-order quasilinear term in \eq{F-evol-reg} is
(even uniformly) monotone and \MMM its \TTT  limit passage (exploiting compact
embedding and Minty's trick or strong convergence) 
is easy, cf.\ e.g.\ \cite[Ch.8]{Roub13NPDE}. \EEE
The limit is a weak solution to the initial-boundary value problem for
\eq{F-evol-reg},  which we denote by $\FF_\varepsilon\in
 H^1(I;  L^2(\varOmega;\R^{d\times d}))\cap  L^\infty 
(I;W^{1,r}(\varOmega;\R^{d\times d}))$.  As this solution is
unique, no extraction of subsequences is actually needed and 
the whole sequence
$\{\FF_k\}_{k\in\N}^{}$ converges to  $\FF_\varepsilon$\TTT; the uniqueness
for $\vv$ given is easy by the uniform monotonicity of the
quasilinear term \MMM and by \TTT  handling the lower-order terms by
Green formula as in \eq{test-FF} \MMM and \TTT  the Gronwall inequality. \EEE

Recalling now that 
\begin{align}
\bigg\|\pdt{\FF_\varepsilon}+(\vv{\cdot}\nabla)\FF_\varepsilon-(\nabla\vv)\FF_\varepsilon
\bigg\|_{L^2( I{\times}\varOmega ;\R^{d\times d})}\le 
  C{\rm e}^{1/(r\varepsilon)}  \,,\end{align}
by comparison  in \eqref{F-evol-reg}  we also obtain
\begin{align}\label{Euler-quasistatic-est3-1}
&\| \varepsilon  {\rm div}(|\Nabla\FF_\varepsilon|^{r-2}\Nabla\FF_\varepsilon^{})\|_{L^2( I{\times}\varOmega ;\R^{d\times d})}^{}\le
C
{\rm e}^{1/(r\varepsilon)}\,.\end{align}
Note that this estimate degenerates  for  
$\varepsilon\to0$.  Still, we have 
that ${\rm div}(|\Nabla\FF_\varepsilon|^{r-2}\Nabla\FF_\varepsilon)\in
L^2( I{\times}\varOmega ;\R^{d\times d})$,  so that  
\eq{F-evol-reg} is solved almost everywhere. In particular, 
  we can legitimately
 test it by   ${\rm div}(|\Nabla\FF_\varepsilon|^{r-2}\Nabla\FF_\varepsilon)$.
Since 
$p>d$,
we have $
p^{-1}+(r^*)^{-1}+(r')^{-1}\le1$,
and thus by the H\"older and Young inequalities, we can estimate
\begin{align}\label{test-Delta-r}
\frac{\d}{\d t}
\int_\varOmega\frac1r|\nabla\FF_\varepsilon|^r\,\d \xx&\le\frac{\d}{\d t}
\int_\varOmega\frac1r|\nabla\FF_\varepsilon|^r\,\d \xx
+
\varepsilon\int_\varOmega|{\rm div}(|\Nabla\FF_\varepsilon|^{r-2}\Nabla\FF_\varepsilon)|^2\,\d \xx
\\[-.0em]\nonumber&
=  -  
\int_\varOmega\nabla\big((\vv{\cdot}\nabla)\FF_\varepsilon-(\nabla\vv)\FF_\varepsilon
  \big)\Vdots\big(|\nabla\FF_\varepsilon|^{r-2}\nabla\FF_\varepsilon\big)\,\d \xx
\\[-.1em]\nonumber&=  - 
\int_\varOmega|\nabla\FF_\varepsilon|^{r-2}(\nabla\FF_\varepsilon{\otimes}\nabla\FF_\varepsilon){:}\ee(\vv)-\frac1r|\nabla\FF_\varepsilon|^r{\rm div}\,\vv
\\[-.4em]&\hspace*{2em}\nonumber
-\big((\nabla\vv)\nabla\FF_\varepsilon+(\nabla^2\vv)\FF_\varepsilon
\big)\Vdots\big(|\nabla\FF_\varepsilon|^{r-2}\nabla\FF_\varepsilon\big)
\,\d \xx
\\&\nonumber
\le C_r
\|\nabla\vv\|_{L^\infty(\varOmega;\R^{d\times d})}^{}
\|\nabla\FF_\varepsilon\|_{L^r(\varOmega;\R^{d\times d\times d})}^r\!
\\[.1em]&\nonumber\hspace*{2em}+C_r
\|\nabla^2\vv\|_{L^p(\varOmega;\R^{d\times d\times d})}^{}
\|\FF_\varepsilon\|_{L^{r^*}(\varOmega;\R^{d\times d})}^{}
\|\nabla\FF_\varepsilon\|_{L^{r}(\varOmega;\R^{d\times d\times d})}^{r-1}
\\[.1em]&\le\nonumber
C_r
\|\nabla\vv\|_{L^\infty(\varOmega;\R^{d\times d})}^{}
\|\nabla\FF_\varepsilon\|_{L^r(\varOmega;\R^{d\times d\times d})}^r
\\[.1em]&\nonumber\hspace*{2em}+
C_rN
\|\nabla^2\vv\|_{L^p(\varOmega;\R^{d\times d\times d})}
\|\FF_\varepsilon\|_{L^2(\varOmega;\R^{d\times d})}
\big(1{+}\|\nabla\FF_\varepsilon\|_{L^{r}(\varOmega;\R^{d\times d\times d})}^r\big)
\\[.1em]&\hspace*{2em}+
C_rN
\|\nabla^2\vv\|_{L^p(\varOmega;\R^{d\times d\times d})}^{}
\|\nabla\FF_\varepsilon\|_{L^{r}(\varOmega;\R^{d\times d\times d})}^r\,,
\nonumber\end{align} 
where we used $p
>d$    in order to get  $\nabla\vv \in W^{1,p}(\varOmega; \R^{d\times d})
\subset  L^\infty(\varOmega;\R^{d\times d})$, as well as the
computation  
\begin{align}\nonumber
&\int_\varOmega\nabla\big((\vv{\cdot}\nabla)\FF_\varepsilon
  \big){:}|\nabla\FF_\varepsilon|^{r-2}\nabla\FF_\varepsilon\,\d \xx
 \\[-.5em]&\hspace{2em}\nonumber=\int_\varOmega|\nabla\FF_\varepsilon|^{r-2}(\nabla\FF_\varepsilon{\otimes}\nabla\FF_\varepsilon){:}\ee(\vv)  
+(\vv{\cdot}\nabla)\nabla\FF_\varepsilon\Vdots|\nabla\FF_\varepsilon|^{r-2}\nabla\FF_\varepsilon\,\d \xx
\\[-.1em]&\hspace{2em}\nonumber
=\int_\varGamma|\nabla\FF_\varepsilon|^r\vv{\cdot}\nn\,d S
+\int_\varOmega\Big(|\nabla\FF_\varepsilon|^{r-2}(\nabla\FF_\varepsilon{\otimes}\nabla\FF_\varepsilon){:}\ee(\vv)
\\[-.8em]&\hspace{14em}\nonumber
-({\rm div}\,\vv)|\nabla\FF_\varepsilon|^r-(r{-}1)|\nabla\FF_\varepsilon|^{r-2}\nabla\FF_\varepsilon\Vdots
(\vv{\cdot}\nabla)\nabla\FF_\varepsilon\Big)\,\d \xx
\\[-.4em]&\hspace{2em}\nonumber
=\int_\varGamma\frac{|\nabla\FF_\varepsilon|^r\!\!}r\ \vv{\cdot}\nn\,d S
+\int_\varOmega|\nabla\FF_\varepsilon|^{r-2}(\nabla\FF_\varepsilon{\otimes}\nabla\FF_\varepsilon){:}\ee(\vv)-({\rm div}\,\vv)\frac{|\nabla\FF_\varepsilon|^r\!\!}r\ \d \xx\,.
\end{align}
Again, the boundary integral  above vanishes since
 
$\vv{\cdot}\nn=0$. For the last inequality in \eq{test-Delta-r}, we
 have used  $\|\FF_\varepsilon\|_{L^{r^*}(\varOmega;\R^{d\times d})}^{}\le
N(\|\FF_\varepsilon\|_{L^2(\varOmega;\R^{d\times d})}^{}+\|\Nabla\FF_\varepsilon\|_{L^{r}(\varOmega;\R^{d\times d})}^{})$,
where $N$ is the norm of the embedding $W^{1,r}(\varOmega)\subset
L^{r*}(\varOmega)$ if $W^{1,r}(\varOmega)$ is endowed with the norm
$\|\cdot\|_{L^2(\varOmega)}+\mbox{$\|\nabla\cdot\|_{L^r(\varOmega;\R^d)}$}$.

 One can thus  apply the Gronwall inequality to \eq{test-Delta-r}. 
 Correspondingly,  
by  using
the 
former estimate in 
\eqref{Euler-quasistatic-est1-2}
and  the regularity of the initial datum  $\Fes_0\in
W^{1,r}(\varOmega;\R^{d\times d})$ 
one obtains the estimates 
\begin{subequations}\label{Euler-quasistatic-est3}
\begin{align}\label{Euler-quasistatic-est3-1+}
&\|\Nabla\FF_\varepsilon\|_{L^\infty(I;L^r(\varOmega;\R^{d\times d\times d}))}\le
     C\ \ \text{ and }
\\&\label{Euler-quasistatic-est3-2}
\|{\rm div}(|\Nabla\FF_\varepsilon|^{r-2}\Nabla\FF_\varepsilon^{})\|_{L^2( I{\times}\varOmega ;\R^{d\times d})}\le C
\varepsilon^{-1/2}\,.
\end{align}\end{subequations} 
The limit passage for
$\varepsilon\to0$ in  the  linear terms is then  straightforward
 and the quasilinear regularizing term in \eq{F-evol-reg}
 converges to $0$  as $\mathscr{O}(\varepsilon^{1/2})$ for
$\varepsilon\to0$ due to \eq{Euler-quasistatic-est3-2}. Alternatively,  one
can observe that  when tested by $\widetilde\SS \UUU \in L^r(I;W^{1,r}(
            \varOmega;\R^{d\times d}))$ \EEE and  by  using
\eq{Euler-quasistatic-est3-1+},  quasilinear regularizing term converges
to $0$  even faster as
$$
\bigg|\int_0^T\!\!\!\int_\varOmega\varepsilon|\Nabla\FF_\varepsilon|^{r-2}\Nabla\FF_\varepsilon\Vdots
\Nabla\widetilde\SS\,\d \xx\d t\bigg|\le
\varepsilon\|\Nabla\FF_\varepsilon\|_{L^r(I{\times}\varOmega;\R^{d\times d\times d})}^{r-1}
\|\Nabla\widetilde\SS\|_{L^r( I{\times}\varOmega;\R^{d\times d\times d})}^{}=\mathscr{O}(\varepsilon)\,.
$$
In any case, the limit for $\varepsilon\to0$ solves the original 
nonregularized  initial-boundary
value problem for \eq{Euler2-diff}. As this equation is linear,
 the  solution is unique  and no extraction of
subsequences is needed in the $\varepsilon \to 0$ limit passage. 

 Estimate \eq{Euler-quasistatic-est2} 
on $\pdt{}\FF_\varepsilon$  does not pass to the limit as
$\varepsilon \to 0$. Still, we can argue by comparison in 
$\pdt{}\Fes=(\Nabla\vv)\Fes
-(\vv{\cdot}\Nabla)\Fes$ 
 and get the estimate  
\begin{align}
\Big\|\pdt{\Fes}\Big\|_{L^1(I;L^r(\varOmega;\R^{d\times d}))}\le C\,.\label{eq:Ft}
\end{align}
In particular, \eq{Euler2-diff} holds a.e.\ on $ I{\times}\varOmega$.
By the embedding $$L^\infty(I;W^{1,r}(\varOmega;\R^{d\times d}))
\cap W^{1,1}(I;L^r(\varOmega;\R^{d\times d}))\subset C(I{\times}\barOmega;\R^{d\times d}),$$
we  also have that  $\Fes\in C(I{\times}\barOmega;\R^{d\times d})$.

The (weak,weak*)-continuity of the mapping $\vv\mapsto\Fes$  is
easy to obtain. Let $\vv_n \to \vv$ weakly in
$L^1(I;W^{2,p}(\varOmega;\R^{d\times d}))$ and let $\FF_n$ be the
corresponding unique solutions of \eqref{Euler2-diff}. Starting from the bound \eqref{eq:Ft} on 
$\FF_n$ in $W^{1,1}(I;L^r(\varOmega;\R^{d\times d}))$ (which indeed
depends on $\| \vv_n
\|_{L^1(I;W^{2,p}(\varOmega;\R^{d\times d}))}$), one applies 
the Aubin-Lions theorem
 obtaining  strong convergence of  $\FF_n$ 
in $L^{1/\epsilon}(I;L^{r^*-\epsilon}(\varOmega;\R^{d\times d}))$ for any $0<\epsilon\le1$.
Then,  we simply pass to the limit in \eq{Euler2-diff} 
in its weak formulation
\eq{def-Euler-weak2}  as $n\to\infty$. 

 Recall that $\Fes_0 \in W^{1,r}(\varOmega;\R^{d\times d})\subset
L^\infty (\varOmega;\R^{d\times d})$.  If   $\det\Fes_0>0$ on $\barOmega$,
$\FF_0^{-1}$ exists and is bounded on $\barOmega$.
 In fact we have that   $\FF_0^{-1}\in
W^{1,r}(\varOmega;\R^{d\times d})$  as 
$$
\nabla\FF_0^{-1}=\nabla\bigg(\frac{{\rm Cof}\FF_0}{\det\FF_0}\bigg)=
\bigg(\frac{{\rm Cof\,}'(\FF_0)}{\det\FF_0}
-\frac{{\rm Cof}(\FF_0){\rm Cof}(\FF_0)}{\det\FF_0^2}\bigg)\nabla\FF_0\in
L^r(\varOmega;\R^{d\times d\times d})\,.
$$
We can then apply the above arguments to  the flow  equation \eq{eq:FF2}
for the  inverse  $\FF^{-1}$,  as well.  
In particular, we obtain  that  $\FF^{-1}$  is  bounded on $ I{\times}\varOmega $,
so that $1/\det\FF$ stays positive and  bounded  away from 0.
\end{proof}

A scalar-valued variant of Lemma~\ref{lem-transport-F+}  holds for
 the continuity equation 
\eq{eq:model1}. Its weak formulation  corresponds to the  
integral identity
\begin{align}
\int_0^T\!\!\!\int_\varOmega\varrho\pdt{v}+\varrho\vv{\cdot}\nabla v\,\d \xx\d t
=\int_\varOmega\varrho_0v(0)\,\d \xx\label{eq:conti}
\end{align}
for any $v\in C^1( I{\times}\barOmega )$ with $v(T)=0$. 
We state this Lemma without proof, for the sake of completeness and later reference. 
 
\begin{lemma}[Flow of $\varrho$] 
\label{lem2}
Let $p>d$  and  $r>2$.
Then, for any $\vv\in L^1(I;W^{2,p}(\varOmega;\R^d))$ with 
$\vv{\cdot}\nn=0$  
 and any $\varrho_0\in W^{1,r}(\varOmega)$, there exists  a unique
  weak
solution $\varrho\in C_{\rm w}(I;W^{1,r}(\varOmega))\cap W^{1,1}(I;L^r(\varOmega))$
to 
\eqref{eq:model1}  in the sense of \eqref{eq:conti}  
and the estimate
\begin{align}\label{rho-evol-est}
\|\varrho\|_{L^\infty(I;W^{1,r}(\varOmega))\,\cap\, W^{1,1}(I;L^r(\varOmega))}^{}\le
\mathfrak{C}\Big(\|\Nabla\vv\|_{L^1(I;W^{1,p}(\varOmega;\R^{d\times d}))}^{}\,,\,
\|\varrho_0\|_{W^{1,r}(\varOmega)}^{}\Big)
\end{align}
holds with some $\mathfrak{C}\in C(\R^2)$. 
Moreover, $\varrho\in C( I{\times}\barOmega )$ and the mapping
\begin{align}
\vv\mapsto\varrho:L^1(I;W^{2,p}(\varOmega;\R^d))\to
L^\infty(I;W^{1,r}(\varOmega))
\end{align}
is (weak,weak*)-continuous.
If in addition $\varrho_0>0$ on $\barOmega$, then $\varrho>0$ on
$ I{\times}\barOmega $ uniformly  with respect to bounded
velocity fields $\vv$, namely,  
for any  $R>0$ there exists  $\delta>0$ such that
\begin{align}
\|\Nabla\vv\|_{L^1(I;W^{1,p}(\varOmega;\R^{d\times d}))}^{}\le R\ \ \Rightarrow\ \
\min_{ I{\times}\barOmega }\varrho\ge\delta\,.
\label{rho>0}\end{align}
\end{lemma}

 We are now in the position of stating the main result of this section. 

\begin{proposition}[Existence and regularity of weak solutions]
\label{prop-Euler-plast}
 Under  assumptions \emph{\eqref{ass}} there exits a weak
solution  $(\varrho,\vv, \bm{F},z,\mu)$  to the
initial-boundary-value problem \emph{\eq{Euler-large-CH-plast}}--\emph{\eq{Euler-BC}} with
\emph{\eq{IC}}  in the sense of  Definition~\emph{\ref{def}}. Moreover,
\begin{align*}
&\FF\in H^1(I;L^r(\varOmega;\R^{d\times d})),\\
&z\in L^2(I;H^1(\varOmega))\cap H^1(I;H^1(\varOmega)^*),\\ 
&\varrho=\rhoR/\det\FF\in H^1(I;L^r(\varOmega)), \\
&\varrho\vv\in L^4(I;W^{1,4}(\varOmega;\R^d))\cap
W^{1,p'}(I;W^{2,p}(\varOmega;\R^d)^*).
\end{align*}
 Eventually,   
the energy dissipation balance
\emph{\eq{Euler-swelling-large-energy}} holds \TTT integrated \EEE
on the time interval $[0,t]$ for any $t\in I$.
\end{proposition}

 The remainder of the paper is devoted to a proof of the latter
existence statement. This hinges upon a nested
regularization and Galerkin space-approximation procedure. In
particular, nonlinearities firstly are replaced by
regularizations. Then, the PDE problem is reduced to an ODE system by
resorting to finite dimensional subspaces. The crucial point here is
that the (weak formulations of the) momentum equation
\eqref{Euler1-large-diff} and of the diffusion equation
\eqref{Euler4-large-diff} will be space discretized. The
continuity equation \eqref{eq:conti} and the flow equation for
$\bm{F}$ will not be space discretized, in order to take advantage
of Lemmas \ref{lem-transport-F+} and \ref{lem2}. 

\begin{proof}
 As mentioned, the proof relies on subsequent approximations and
is here \UUU divided \EEE  into steps, for better clarity.  

\medskip\noindent{\it Step 1: Regularization}.
\UUU Since $r>d$, we can choose $\varepsilon>0$  small enough so that
all fields $\FF$ fulfilling the formal estimate  \eq{est+Fes} satisfy
\EEE 
\begin{align}
&
\det\FF>\varepsilon\ \ \ \ \text{ and }\ \ \ \ |\FF|<\frac1\varepsilon
\ \ \text{ a.e.\ on }\  I{\times}\varOmega \,.
\label{Euler-quasistatic-est-formal4}
\end{align}
\UUU Correspondingly, we may perform a regularization of \EEE 
the stress in 
\eq{Euler-large-CH-plast} by considering a smooth cut-off
$\varphi_\varepsilon(\cdot,z)$ of the original stored  \OK energy \EEE
density   $\varphi$ defined as
\begin{align}\label{cut-off-general}
&\varphi_\varepsilon(\Fe,z)=\chi_\varepsilon(\lambda(z)\Fe)\varphi(\Fe,z)
\\
&\text{ with }\ \chi_\varepsilon(\FF)
=\begin{cases}
\qquad\qquad1&\hspace{-8em}\text{for $\det \FF\ge\varepsilon$ and $|\FF|\le1/\varepsilon$,}
\\
\qquad\qquad0&\hspace{-8em}\text{for $\det \FF\le\varepsilon/2$ or $|\FF|\ge2/\varepsilon$,}
\\
\displaystyle{\Big(\frac{3}{\varepsilon^2}\big(2\det\FF-\varepsilon\big)^2
-\frac{2}{\varepsilon^3}\big(2\det\FF-\varepsilon\big)^3\Big)\,\times}\!\!&
\\[.2em]
\qquad\qquad\displaystyle{\times\,\big(3(\varepsilon|\FF|-1)^2
-2(\varepsilon|\FF|-1)^3\big)}\!\!&\text{otherwise}.
\end{cases}
\nonumber\end{align}
 We moreover make use of the notation  
$\widehat\varphi_\varepsilon(\FF,z)=\chi_\varepsilon(\FF)
\widehat\varphi(\FF,z)$.
Note that also
$\varphi_\varepsilon,\widehat\varphi_\varepsilon\in C^1(\R^{d\times d}\times\R)$ if
$\varphi\in C^1(\R^{d\times d}\times\R)$. \UUU Moreover, \EEE 
$[\widehat\varphi_\varepsilon]_\FF'$, the
Cauchy stress $(\FF,z)\mapsto\TT_\varepsilon=
[\widehat\varphi_\varepsilon]_\FF'(\FF,z)\FF^\top+\widehat\varphi_\varepsilon(\FF,z)\bbI$,
\UUU and \EEE the driving pressure $\pi_\varepsilon=[\widehat\varphi_\varepsilon]_z'$
are bounded, continuous. \UUU In fact, $\TT_\varepsilon$ and
$\pi_\varepsilon$ vanish as an effect of the choice of
$\chi_\varepsilon$ \TTT if $\FF$ ``substantially'' violates the bounds
\eq{Euler-quasistatic-est-formal4}, specifically \UUU if 
$\det\FF\le\varepsilon/2$ or $|\FF|\geq
2/\varepsilon$. \EEE  
It is also important to  notice  that the strong convexity of
$\widehat\varphi(\FF,\cdot)$ is  
{\it not}  inherited  by $\widehat\varphi_\varepsilon(\FF,\cdot)$, which is why
we 
 are forced to resort to a  regularization of the 
diffusion equation  
\eq{Euler4-large-CH-reg} below. 

 The multivalued mapping  $N_{[0,1]}^{}(\cdot)$ in
\eq{Euler4-large-diff} is approximated via  
the standard Yosida approximation
\begin{align}
\mathcal{N}_k(z)=\begin{cases}
k(z{-}1)&\text{ if }z>1,\\[-.2em]
\ \ \ 0&\text{ if }0\le z\le1,\\[-.2em]
\ \ kz&\text{ if }z<0.
\end{cases}
\end{align}
 Note that $k \in \N$ is the \UUU index \EEE  of the Galerkin approximation
of the momentum equation as well, see Step2 below. 

We moreover regularize the singular nonlinearity
$1/\det(\cdot)$,  showing up  
in the right-hand-side  of  the momentum
equation, although simultaneously the mass-density continuity equation is
considered for the inertial term. To this aim, 
we introduce the short-hand notation
\begin{align}\label{cut-off-det}
\mbox{$\det_\varepsilon$}\FF:=\max(\det\FF,\varepsilon)\,,
\end{align}
 Eventually, we  regularize also the diffusion equation for
$z$.  
Altogether, the regularized system  reads as follows 
\begin{subequations}\label{Euler-large-CH-reg}
\begin{align}
\label{Euler0=hypoplast}&
\pdt\varrho=-{\rm div}(\varrho\vv)
\,,
\\\label{Euler1-large-CH-reg}
     &\pdt{}(\varrho\vv)+{\rm div}(\varrho\vv{\otimes}\vv)
  ={\rm div}(\TT_\varepsilon{+}
  \Sv)
 +\frac{\rhoR\GRAVITY}{\det_\varepsilon\Fes}
     \ \ \
     \\&\hspace*{7.2em}\nonumber
\text{with }\ \,\TT_\varepsilon=
[\widehat\varphi_\varepsilon]_\FF'(\FF,z)\FF^\top+\widehat\varphi_\varepsilon(\FF,z)\bbI
\nonumber\\&\hspace*{7.2em}
\text{and }\ \ \Sv=\zeta_\ee'(z;\EE(\vv))-{\rm div}
\big(\nu|\nabla\EE(\vv)|^{p-2}\nabla\EE(\vv)\big)\,,
\nonumber\\\label{Euler2-large-CH-reg}
&\DT\Fes=(\Nabla\vv)\Fes
\,,
\\&
\DT z={\rm div}\Big(\widehat\M(\FF,z)\nabla\mu
+\big(1{-}\chi_\varepsilon(\FF)\big)\nabla z\Big)
\ \ \
\text{ with }\ \ 
\mu=[\widehat\varphi_\varepsilon]_z'(\FF,z)
+\mathcal{N}_k(z)\,.
\label{Euler4-large-CH-reg}\end{align}\end{subequations}
The boundary conditions for \eqref{Euler1-large-CH-reg} are as in
\eqref{Euler-BC} while the condition for the diffusion equation
(i.e., the last condition in \eqref{Euler-BC}) is now modified as
\begin{align}\label{BC-regular}
\Big(\widehat\M(\FF,z)\nabla\mu
+\big(1{-}\chi_\varepsilon(\FF)\big)\nabla z\Big){\cdot}\nn
+\varkappa\mu
+\big(1{-}\chi_\varepsilon(\FF)\big)z=h\,.
\end{align} 
 Note that the  terms with factor $1{-}\chi_\varepsilon(\FF)$
in
\eqref{Euler4-large-CH-reg} and \eqref{BC-regular}  vanish 
if $\FF$ complies with the bounds
\eqref{Euler-quasistatic-est-formal4}.  On   the other
hand,  they ensure the  strong monotonicity of the diffusion operators
and the coercivity of the boundary conditions, even when the approximate solution
violates these bounds and thus the cut-off $[\widehat\varphi_\varepsilon]_z'$
may degenerate.

\medskip\noindent{\it Step 2:  Galerkin approximation}.  
 We perform a   
Galerkin approximation separately of  the momentum equation \eq{Euler1-large-CH-reg}
and of  the diffusion equation for $z$ \eq{Euler4-large-CH-reg}. 
 On the other hand, we do not approximate in space  the  continuity
 equation  
\eq{Euler0=hypoplast}
and  the flow equation \eq{Euler2-large-CH-reg} for $\bm{F}$  
 but rather rely respectively on  
 Lemmas~\ref{lem2} and \ref{lem-transport-F+} for their
weak solutions.  The Galerkin approximations of equations
\eq{Euler1-large-CH-reg} and \eq{Euler4-large-CH-reg} are
kept independent, in order to be \TTT able to pass \EEE
separately to the limit in Steps~6 and 4, respectively. 

Specifically, we use a nested finite-dimensional subspaces $\{V_k\}_{k\in\N}$
whose union is dense in $W^{2,p}(\varOmega;\R^d)$ for  the momentum equation 
\eq{Euler1-large-CH-reg}.  Note that
these spaces  are indexed by the same $k\in\N$ used
in \eq{Euler4-large-CH-reg} for the regularization of the
normal-cone mapping.  In addition, we perform  a Galerkin
approximation of 
the  diffusion  equation \eq{Euler4-large-CH-reg} by
using  a second collection of  nested
finite-dimensional subspaces $\{Z_l\}_{l\in\N}$  whose union is dense in
$H^1(\varOmega)$.  
Without loss of generality, we may assume $\vv_0\in V_1$
and $z_0\in Z_1$.

The  space approximation of the   
solution of the regularized system
\eq{Euler-large-CH-reg} will be denoted by
$$(\varrho_{kl},\vv_{kl},\FF_{kl},z_{kl})
:I\to W^{1,r}(\varOmega)\times V_k\times W^{1,r}(\varOmega;\R^{d\times d})
\times Z_l.$$
Existence of  such space-approximated  solution  can be
obtained via the standard existence theory for first-order systems of 
ordinary differential equations:  local-in-time existence follows
from smoothness, also in connection with  
Lemmas~\ref{lem-transport-F+} and \ref{lem2}.  Then, global
existence on the whole time interval $[0,T]$ results from the standard
successive-prolongation argument, on the basis of the uniform-in-time
estimates proved below.

Let us once again stress that the continuity equation
\eq{Euler0=hypoplast} is not space discretized. This allows us to test
it by   
$|\vv_{kl}|^2$  
 so that  
identity \eq{inertial}  
is at disposal also
 at  the Galerkin level.  
On the other hand, it is to be
emphasized that also the equation for $\mu$ in \eq{Euler4-large-CH-reg}
is not  space  discretized:  the corresponding 
$\mu_{kl}$ is  therefore  not valued in $Z_l$ and  thus is not a legitimate test function for the diffusion
equation \eq{Euler4-large-CH-reg}.

\medskip\noindent{\it Step 3: First a-priori estimates}.  A basic
estimate follows by testing the  Galerkin approximation of the
momentum equation  
\eq{Euler1-large-CH-reg} by $\vv_{kl}$,  taking advantage of
the  (not  discretized!) continuity
equation \eq{Euler0=hypoplast} tested by $|\vv_{kl}|^2/2$,
and  by testing the  Galerkin approximation of the diffusion
equation  \eq{Euler4-large-CH-reg} by  $z_{kl}$.

 The continuity equation  \eq{Euler0=hypoplast} tested   by $|\vv_{kl}|^2/2$ can be used
 in  \eq{inertial},  here   written in terms of
$\varrho_{kl}$ and $\vv_{kl}$, 
in order to  exploit the formulas  \eq{rate-of-kinetic}--\eq{calculus-convective-in-F}
to obtain the rate of kinetic energy.  A crucial observation is
  that, due to the  presence of the  
cut-offs $\det_\varepsilon$ and $\varphi_\varepsilon$, the equations
(\ref{Euler-large-CH-reg}a--c) can be estimated independently of $z$,
i.e., independently  from the estimate of the diffusion equation  \eq{Euler4-large-CH-reg}. Specifically,  from
 the Galerkin approximation of  \eq{Euler1-large-CH-reg}
tested by $\vv_{kl}$  we obtain the identity
\begin{align}\label{Euler-swelling-large-energy++}
  &\frac{\d}{\d t}
  \int_\varOmega\!
\frac{\varrho_{kl}}2
|\vv_{kl}|^2
\,\d \xx
+\int_\varOmega 
\zeta_\ee'(z;\EE(\vv_{kl})){:}\EE(\vv_{kl})
+\nu|\nabla\EE(\vv)|^p
\,\d \xx
\\[-.5em]&\hspace{9.5em}
=\int_\varOmega
\frac{\rhoR\GRAVITY}{\det_\varepsilon\Fes_{kl}\!\!}
{\cdot}\vv_{kl}-
\TT_{\varepsilon,kl}
{:}\EE(\vv_{kl})
+\int_\varGamma\TRACTION{\cdot}\vv_{kl}
\,\d S
\nonumber\end{align}
where 
\begin{align}
&\nonumber\\[-2.5em]
&\TT_{\varepsilon,kl}=
[\widehat\varphi_\varepsilon]_{\FF}'(\FF_{kl},z_{kl})\FF_{kl}^\top
+\widehat\varphi_\varepsilon(\FF_{kl},z_{kl})\bbI
\,.
\end{align}

 Due to  Lemmas~\ref{lem-transport-F+} and \ref{lem2} with $\vv=\vv_{kl}$
and with the fixed initial conditions $\FF_0$ and $ \varrho_0$, 
 we may  define the nonlinear operators
$\mathfrak{F}:I\times L^p(I;W^{2,p}(\varOmega;\R^d))
\to W^{1,r}(\varOmega;\R^{\d\times d})$ and
$\mathfrak{R}:I\times L^p(I;W^{2,p}(\varOmega;\R^d))
\to W^{1,r}(\varOmega)$ by 
\begin{align}
\FF_{kl}(t)=\mathfrak{F}\big(t,\vv_{kl}
\big)\ \ \text{ and }\ \
\varrho_{kl}(t)=\mathfrak{R}\big(t,\vv_{kl}
\big)\,.
\end{align} 

Since  we have that  $p\ge2$, we can estimate
\begin{align}\label{est-fv}
&\int_\varGamma\!\TRACTION{\cdot}\vv_{kl}\,\d S
                         \le\|\TRACTION\|_{L^1(\varGamma;\R^d)}\|\vv_{kl}\|_{L^\infty(\varGamma;\R^d)}
  \\&\nonumber
\quad \le N\|\TRACTION\|_{L^1(\varGamma;\R^d)}\big(\|\vv_{kl}\|_{L^2(\varGamma;\R^d)}^{}\!+
\|\nabla\EE(\vv_{kl})\|_{L^p(\varOmega;\R^{d\times d\times d})}^{}\big)
\\&\nonumber \quad  \le C\|\TRACTION\|_{L^1(\varGamma;\R^d)}^{p'}
+\|\TRACTION\|_{L^1(\varGamma;\R^d)}\big(1+\|\vv_{kl}\|_{L^2(\varOmega;\R^d)}^2\big)
+\delta\|\nabla\EE(\vv_{kl})\|_{L^p(\varOmega;\R^{d\times d\times d})}^p
\\& \quad \le C\|\TRACTION\|_{L^1(\varGamma;\R^d)}^{p'}\!
+\|\TRACTION\|_{L^1(\varGamma;\R^d)}\bigg(1{+}
\frac{\!\|\sqrt{\varrho_{kl}}\vv_{kl}\|_{L^2(\varOmega;\R^d)}^2\!}{\sqrt{\min\varrho_{kl}}}\bigg)
+\delta\|\nabla\EE(\vv_{kl})\|_{L^p(\varOmega;\R^{d\times d\times d})}^p\,,
\nonumber\end{align}
where $N$ depends on  the norm of the trace operator
$W^{2,p}(\varOmega)\to L^\infty(\varGamma)$ and  the
Korn-inequality  constant,
while $C$ depends on $N$ and $\delta>0$, which can be chosen arbitrarily.

By the 
Gronwall inequality, we obtain the estimates 
\begin{subequations}\label{Euler-quasistatic-est1}
\begin{align}\label{Euler-quasistatic-est1-1}
&\|\EE(\vv_{kl})\|_{
L^2(I;W^{1,p}(\varOmega;\R^{d\times d}))}^{}\le C\,\ \ \text{ and }\ \
\big\|\sqrt{\varrho_{kl}}\vv_{kl}\big\|_{L^\infty(I;L^2(\varOmega;\R^d))}\le C\,,
\intertext{and, since $\varrho_{kl}$  is uniformly bounded away
  from $0$,  from   \eqref{rho>0} together with
\eq{Euler-quasistatic-est1-1}, 
we also have  that } 
&
\|\vv_{kl}\|_{L^\infty(I;L^2(\varOmega;\R^d))}\le C\,.
\end{align}\end{subequations}

Next, we use the strong convexity of
$\widehat\varphi(\FF,\cdot)$, cf.\ \eq{ass-phi},
in order to  drop momentarily  
equation
$\mu_{kl}=[\widehat\varphi_\varepsilon]_z'(\FF_{kl},z_{kl})+\UUU \mathcal{N}_k\EEE(z_{kl})$
 which holds  a.e.\ on $ I{\times}\varOmega. $ 
Indeed, this equation  
should otherwise be tested by
$\DT z_{kl}$, which
would not be a legitimate test  at the Galerkin-approximation 
level. 
 By computing the gradient , we have 
\begin{align}
\nabla\mu_{kl}
&=
\big([\widehat\varphi_\varepsilon]_{zz}''(\FF_{kl},z_{kl})+\xi_{kl}\big)\nabla z_{kl}
+[\widehat\varphi_\varepsilon]_{\FF z}''(\FF_{kl},z_{kl})\nabla\FF_{kl}
\ \ \text{ with }\xi_{kl}\in
\UUU \mathcal{N}_k'\EEE(z_{kl})\,.
\end{align} 
Note that $\mathcal{N}_k\in W^{2,\infty}(\R)$ and the (generalized)
\UUU derivative $\mathcal{N}_k'$ \EEE indeed jumps (i.e., is set-valued) at $z=0$
and $z=1$.  On the other hand, we nevertheless have that    $0\le\xi_{kl}\le1/k$.
Substituting  this   into \eq{Euler4-large-CH-reg}, we obtain an
initial-boundary-value problem for $z_{kl}$,  namely, (the
Galerkin \UUU approximation \EEE of)  
\begin{align}\label{evol-of-z}
&\DT z_{kl}={\rm div}\,\jj_{kl}\ \ \ \text{ with}
\\&
\jj_{kl}=\Big(\!\!\!\!\linesunder{
\widehat\M(\FF_{kl},z_{kl})
\big([\widehat\varphi_\varepsilon]_{zz}''(\FF_{kl},z_{kl}){+}\xi_{kl}\big)+
1{-}\chi_\varepsilon(\FF_{kl})}{$=:\mathfrak{m}(\FF_{kl},z_{kl})$,
``uniformly" positive}{with respect to $(\FF_{kl},z_{kl})$}
\!\!\!\!\Big)\nabla z_{kl}
+[\widehat\varphi_\varepsilon]_{\FF,z}''(\FF_{kl},z_{kl})\nabla\FF_{kl}
\nonumber\end{align}
and with the boundary condition
$\jj_{kl}{\cdot}\nn+\varkappa\mu_{kl}=h$. 
It is now  allowed  to test
\eq{evol-of-z} in its Galerkin approximation by $z_{kl}$, which 
leads to the identity
\begin{align}\label{est-of-z}
&\frac12\frac{\d}{\d t}\|{z_{kl}}\|_{L^2(\varOmega)}^2
+\int_\varOmega\mathfrak{m}(\FF_{kl},z_{kl})|\nabla z_{kl}|^2\,\d \xx
+\int_\varGamma\varkappa\mu_{kl}z_{kl}\,\d S
\\&\nonumber=
\int_\varOmega[\widehat\varphi_\varepsilon]_{\FF,z}''(\FF_{kl},z_{kl})
\Vdots\big(\nabla\FF_{kl}{\otimes}\nabla z_{kl}\big)
-(\vv_{kl}{\cdot}\nabla z_{kl})z_{kl}\,\d \xx+\int_\varGamma\varkappa h z_{kl}
+(1{-}\chi_\varepsilon(\FF_{kl}))z_{kl}^2\,\d S
\\&
=\int_\varOmega[\widehat\varphi_\varepsilon]_{\FF,z}''(\FF_{kl},z_{kl})
\Vdots\big(\nabla\FF_{kl}{\otimes}\nabla z_{kl}\big)
+\frac{|z_{kl}|^2\!}2\,{\rm div}\,\vv_{kl}
\,\d \xx+\int_\varGamma\varkappa h z_{kl}-
\frac{|z_{kl}|^2\!}2\,\vv_{kl}{\cdot}\nn\,\d S\,,
\nonumber\end{align}
where also the Green formula in  $\varOmega$ has been
used. It is important that  the term 
$\varkappa\mu_{kl}z_{kl}=
\varkappa\chi_\varepsilon(\FF_{kl})\widehat\varphi_z'(\FF_{kl},z_{kl})z_{kl}
+(1{-}\chi_\varepsilon(\FF_{kl}))z_{kl}^2$ can be estimated from below by
$\delta|z_{kl}|^2-1/\delta$
for sufficiently small $\delta>0$, depending on the strong convexity of
$\widehat\varphi(\FF,\cdot)$, cf.\ \eqref{ass-phi}, so that the boundary
term $\int_\varGamma\varkappa h z_{kl}\,\d S$ in \eqref{est-of-z} can be estimated
by using also the coercive left-hand-side term
$\int_\varGamma\varkappa\mu_{kl}z_{kl}\,\d S$. Using the boundary condition
$\vv_{kl}{\cdot}\nn=0$ and the Gronwall and the H\"older inequalities, we obtain
the estimate
\begin{align}
&\label{Euler-quasistatic-est1-5}
\|z_{kl}\|_{L^\infty(I;L^2(\varOmega))\,\cap\,L^2(I;H^1(\varOmega))}^{}\le C\,.
\end{align}
From this, we also obtain an information about
$\mu_{kl}=[\widehat\varphi_\varepsilon]_z'(\FF_{kl},z_{kl})+\UUU \mathcal{N}_k\EEE(z_{kl})$:
\begin{align}
&\label{Euler-quasistatic-est1-6}
\|\mu_{kl}\|_{L^\infty(I;L^2(\varOmega))\,\cap\,L^2(I;H^1(\varOmega))}^{}\le Ck\,.
\end{align}

\medskip\noindent{\it Step 4: Limit passage  for 
  $l\to\infty$}.  By \UUU
  \TTT the obtained a-priori estimates and the sequential \UUU
  weak* compactness of 
  \TTT balls \UUU in the involved spaces,
  \TTT we can standardly use \EEE
the Banach selection principle \TTT \cite[Chap.\,III,\,Thm.\,3]{Bana32TOL}
(i.e., a special form of the Alaoglu-Bourbaki principle devised later for
nonmetrizable situations) and \EEE 
extract some not relabeled subsequence and
$(\varrho_k,\vv_k,\FF_k,z_k,\mu_k):I\to W^{1,r}(\varOmega)\times
V_k\times W^{1,r}(\varOmega;\R^{d\times d})\times H^1(\varOmega)^2$
such that
\begin{subequations}\label{Euler-weak-sln}
\begin{align}
&\!\!\varrho_{kl}\to\varrho_k&&\text{weakly* in $\
L^\infty(I;W^{1,r}(\varOmega))\,\cap\,W^{1,p}(I;L^r(\varOmega))$}\,,
\\\label{Euler-weak-sln-v}
&\!\!\vv_{kl}\to\vv_k&&\text{weakly* in $\
L^\infty(I;L^2(\varOmega;\R^d))\cap
L^2(I;W^{2,p}(\varOmega;\R^d))$,}\!\!&&
\\
&\!\!\FF_{kl}\to\FF_k
\!\!\!&&\text{weakly* in $\ 
L^\infty(I;W^{1,r}(\varOmega;\R^{d\times d}))\,\cap\,
H^1(I;L^2(\varOmega;\R^{d\times d}))$,
}\!\!
\\
&\!\!z_{kl}\to z_k&&\text{weakly* in $\ L^\infty(I;L^2(\varOmega))\cap L^2(I;H^1(\varOmega))$}\,,
\\
&\!\!\mu_{kl}\to\mu_k&&\text{weakly* in $\ L^\infty(I;L^2(\varOmega))\cap L^2(I;H^1(\varOmega))$}
\,.
\end{align}\end{subequations}
 Recalling that $r>d$, by  the Aubin-Lions  Lemma we also have that 
\begin{subequations}\label{Euler-weak+}
\begin{align}\label{rho-conv}
&\varrho_{kl}\to\varrho_{k} 
\hspace*{2em}\text{strongly in }  C(I{\times}\barOmega) 
\end{align}
and $\FF_{kl}\to\FF_k$ strongly in
$ C(I{\times}\barOmega; \R^{d\times d})$. 
 By comparison in the equation in \eq{Euler4-large-CH-reg} we obtain a
bound on  $\pdt{}z_{kl}$,  implying that 
\begin{align}\label{z-conv}
z_{kl}\to z_k\qquad\text{strongly in $L^s( I{\times}\varOmega )$ for
\TTT any \EEE $1\le s<2+4/d$},
\end{align}
cf.\ \cite[Ch.8]{Roub13NPDE}.
Thus, by the continuity of the corresponding Nemytski\u{\i}
(or here simply superposition) mappings, also the 
conservative part of the regularized Cauchy stress and the 
diffusivity 
and the regularized pore pressure in the diffusion  equation 
converge,  namely, 
\begin{align}
&
\TT_{\varepsilon,kl}\to\TT_{\varepsilon,k}=
[\widehat\varphi_\varepsilon]_{\FF}'(\FF_k,z_k)\FF_k^\top+\widehat\varphi_\varepsilon(\FF_k,z_k)\bbI
\hspace*{-0em}&&\hspace*{-1em}\text{strongly in $L^c( I{\times}\varOmega ;\R^{d\times d})$,}
\\&\label{m-strongly}
\widehat\M(\FF_{kl},z_{kl})\to\widehat\M(\FF_k,z_k)
&&\hspace*{-1em}\text{strongly in $L^c( I{\times}\varOmega )$,}
\\&\label{pi-strongly}
[\widehat\varphi_\varepsilon]_z'(\FF_{kl},z_{kl})
\to[\widehat\varphi_\varepsilon]_z'(\FF_k,z_k)
&&\hspace*{-1em}\text{strongly in $L^c( I{\times}\varOmega )$,}
\end{align}\end{subequations}
for any $1\le c<\infty$. It is important to notice that
\begin{align}\label{nabla-rho.v}
\nabla(\varrho_{kl}\vv_{kl})=
\nabla\varrho_{kl}{\otimes}\vv_{kl}+\varrho_{kl}\nabla\vv_{kl}
\end{align}
is bounded in $L^{ \infty}(I;L^r(\varOmega;\R^{d\times d}))$ due to the already obtained
bounds 
\eq{Euler-quasistatic-est1-2} and \eq{rho-evol-est}. 
Therefore, $\varrho_{kl}\vv_{kl}$ converges  weakly*  in
$L^{ \infty }(I;W^{1,r}(\varOmega;\R^d))$.  In fact,  the limit of
$\varrho_{kl}\vv_{kl}$ can be identified as $\varrho_{k}\vv_k$ because we already
showed that $\varrho_{kl}$ converges strongly in \eq{rho-conv} and $\vv_{kl}$
converges weakly due to \eqref{Euler-weak-sln-v}. 

 By comparison, we also obtain  
some information about
$\pdt{}(\varrho_{kl}\vv_{kl})$.  Note indeed that  
\eq{inertial}  still holds   for the semi-discretized system since
the continuity equation  has not been space-discretized.  Specifically, we have
\begin{align}\label{est-of-DT-rho.v}
\pdt{}(\varrho_{kl}\vv_{kl})&=\varrho_{kl}\DT\vv_{kl}
-{\rm div}(\varrho_{kl}\vv_{kl}{\otimes}\vv_{kl})
\\&=\varrho_{kl}\DT\vv_{kl}-\varrho_{kl}(\vv_{kl}{\cdot}\nabla)\vv_{kl}
-\varrho_{kl}({\rm div}\,\vv_{kl})\vv_{kl}-(\vv_{kl}{\cdot}\nabla\varrho_{kl})\vv_{kl}\,.\nonumber
\end{align}
 We may hence compare in \eq{Euler1-large-CH-reg} in order
to obtain a bound on $\varrho_{kl}\DT\vv_{kl}$.   
By the compact embedding $L^\infty(I;V_k)\cap W^{1,p'}(I;V_k)\subset L^\infty(I;V_k)$,
we have 
\begin{align}\label{rho-v-conv}
&\varrho_{kl}\,\vv_{kl}\to\varrho_{k}\vv_k
&&\hspace*{-1em}\text{strongly in }L^{ c }(
   I{\times}\varOmega ;\R^d)\ \ \text{  for all  $1\le c<4$.}\,
\intertext{Since obviously $\vv_{kl}=(\varrho_{kl}\vv_{kl})(1/\varrho_{kl})$,
thanks to \eq{Euler-quasistatic-est1-2}, \eq{rho-conv}, and
\eq{rho-v-conv}, we  also have that }
&\vv_{kl}\to\vv_k
&&\hspace*{-1em}\text{strongly in }L^c( I{\times}\varOmega ;\R^d)\ \ \text{ with
any $1\le c<4$,} 
\,.
\end{align}

The convergences  
(\ref{Euler-weak+}d,e)  allow to pass  to the limit for
$\l\to\infty$ in the regularized  
diffusion equation \eq{Euler4-large-CH-reg}.
The limit passage in the evolution  equations  \eq{Euler0=hypoplast}
and \eq{Euler2-large-CH-reg}  follows from  Lemmas \ref{lem-transport-F+}
and \ref{lem2}.

For the limit passage in the momentum equation, one uses the monotonicity of
the dissipative stress $\DD$, i.e., the
monotonicity of the quasilinear operator
$$\vv\mapsto{\rm div}({\rm div}(|\Nabla\EE(\cdot)|^{p-2}\Nabla\EE(\vv))
-\zeta_{\ee}'(z,\EE(\vv))),$$  and one employs  weak
convergence,  in combination with  the
so-called Minty trick. We take $\widetilde\vv\in H^1(I;V_k)$ and test the
momentum equation by $\vv_{kl}{-}\widetilde\vv$.
 Note that one has  
\begin{align}\label{conv-of-inirtia.v}
&\int_0^T\!\!\!\int_\varOmega\varrho_{kl}\DT\vv_{kl}{\cdot}\widetilde\vv\,\d \xx\d t
=\int_0^T\!\!\!\int_\varOmega
\Big(\pdt{}(\varrho_{kl}\vv_{kl})+{\rm div}(\varrho_{kl}\vv_{kl}{\otimes}\vv_{kl})
\Big){\cdot}\widetilde\vv\,\d \xx\d t
\\[-.4em]&\ \ \ \ \nonumber=\int_\varOmega\varrho_{kl}(T)\vv_{kl}(T){\cdot}\widetilde\vv(T)
-\varrho_0\vv_0{\cdot}\widetilde\vv(0)\,\d \xx
-\!\int_0^T\!\!\!\int_\varOmega\varrho_{kl}\vv_{kl}{\cdot}\pdt{\widetilde\vv}
+(\varrho_{kl}\vv_{kl}{\otimes}\vv_{kl}){:}\nabla\widetilde\vv\,\d \xx\d t
\\&\ \ \ \ \nonumber\to\int_\varOmega\varrho_k(T)\vv_k(T){\cdot}\widetilde\vv(T)
-\varrho_0\vv_0{\cdot}\widetilde\vv(0)\,\d \xx
-\!\int_0^T\!\!\!\int_\varOmega\varrho_k\vv_k{\cdot}\pdt{\widetilde\vv}
+(\varrho_k\vv_k{\otimes}\vv_k){:}\nabla\widetilde\vv\,\d \xx\d t
\\&\ \ \ \
=\int_0^T\!\!\!\int_\varOmega\varrho_k\DT\vv_k{\cdot}\widetilde\vv\,\d \xx\d t\,.
\nonumber\end{align}
Here, we have used the fact that the term  $\varrho_{kl}(T) $ is
also bounded in 
$  W^{1,r}(\varOmega)$ 
 and  $\vv_{kl}(T)$  is bounded in   $ L^2(\varOmega;\R^d)$, together with some information
about the time derivatives $\pdt{}\varrho_{kl}=-{\rm div}(\varrho_{kl}\vv_{kl})$ and
$\pdt{}(\varrho_{kl}\vv_{kl})$, cf.\ \eq{est-of-DT-rho.v}, so that
we can identify the weak limit of $\varrho_{kl}(T)\vv_{kl}(T)$. 
 We have hence obtained that 
\begin{align}
&&&\varrho_{kl}(T)\vv_{kl}(T)\to\varrho_{k}(T)\vv_k(T)&&\text{weakly in $\ L^2(\varOmega;\R^d)$.}&&&&
\label{Euler-weak-rho.v(T)}\end{align}
In \eqref{conv-of-inirtia.v}, we have relied on \eq{rho-v-conv} and on
 the fact  that
$\pdt{}\widetilde\vv$ is well defined  at   the Galerkin
level and that the continuity equation is not discretized, so that  
the identity \eqref{inertial}  holds  even for the semi-discrete
problem. 
This is to be used  in  the following calculations
\begin{align}\label{strong-hyper+}
&0\le
\limsup_{l\to\infty}\bigg(
\int_0^T\!\!\!\int_\varOmega\!\Big(
\nu\big(|\nabla\EE(\vv_{kl})|^{p-2}\nabla\EE(\vv_{kl})
-|\nabla\EE(\widetilde\vv)|^{p-2}\nabla\EE(\widetilde\vv)\big)\Vdots
  \nabla\EE(\vv_{kl}{-}\widetilde\vv)
  \\[-.2em]&\hspace*{3em}\nonumber
  +\big(\pl_\EE'\zeta(z_{kl},\EE(\vv_{kl})){-}\pl_\EE'\zeta(z_{kl},\EE(\widetilde\vv))\big)
{:}\EE(\vv_{kl}{-}\widetilde\vv)
\Big)\,\d \xx\d t\bigg)
\\[-.2em]&=\nonumber
 \limsup_{l\to\infty}\bigg(
 \int_0^T\!\!\!\int_\varOmega
 \bigg(\varrho_{kl}
 (\GRAVITY{-}\DT\vv_{kl})
 {\cdot}(\vv_{kl}{-}\widetilde\vv)
  -
\TT_{\varepsilon,kl}{:}\nabla(\vv_{kl}{-}\widetilde\vv)
 -\pl_\EE'\zeta(z_{kl},\EE(\widetilde\vv)){:}\EE(\vv_{kl}{-}\widetilde\vv)
 \\[-.3em]&\nonumber \hspace*{3em}
-\nu\big(|\nabla\EE(\widetilde\vv)|^{p-2}\nabla\EE(\widetilde\vv)\big)\Vdots
 \nabla\EE(\vv_{kl}{-}\widetilde\vv)
\Big)\,\d \xx\d t
+\int_0^T\!\!\!\int_\varGamma\TRACTION{\cdot}(\vv_{kl}{-}\widetilde\vv)\,\d S\d t
\bigg)%
\\[-.2em]&=\nonumber
 \limsup_{l\to\infty}\bigg(
 \int_0^T\!\!\!\int_\varOmega
 \bigg(\varrho_{kl}\GRAVITY{\cdot}(\vv_{kl}{-}\widetilde\vv)
 +\varrho_{kl}\DT\vv_{kl}{\cdot}\widetilde\vv
  -
\TT_{\varepsilon,kl}{:}\nabla(\vv_{kl}{-}\widetilde\vv)
 \\[-.2em]&\nonumber\hspace*{3em}
 -\pl_\EE'\zeta(z_{kl},\EE(\widetilde\vv)){:}\EE(\vv_{kl}{-}\widetilde\vv)
-\nu\big(|\nabla\EE(\widetilde\vv)|^{p-2}\nabla\EE(\widetilde\vv)\big)\Vdots
 \nabla\EE(\vv_{kl}{-}\widetilde\vv)
\Big)\,\d \xx\d t
\\[-.1em]&\hspace*{3em}\nonumber
+\int_0^T\!\!\!\int_\varGamma\TRACTION{\cdot}(\vv_{kl}{-}\widetilde\vv)\,\d S\d t
-\int_\varOmega\frac{\varrho_{kl}(T)}2|\vv_{kl}(T)|^2-\frac{\varrho_0}2|\vv_0|^2\,\d \xx
\bigg)%
\\[-.2em]&\le\nonumber
\int_0^T\!\!\!\int_\varOmega
 \bigg(\varrho_k\GRAVITY{\cdot}(\vv_k{-}\widetilde\vv)
 +\varrho_k\DT\vv_k{\cdot}\widetilde\vv
  -
\TT_{\varepsilon,k}{:}\nabla(\vv_k{-}\widetilde\vv)
 -\pl_\EE'\zeta(z_k,\EE(\widetilde\vv)){:}\EE(\vv_k{-}\widetilde\vv)
 \\[-.3em]&\nonumber\hspace*{3em}
-\nu\big(|\nabla\EE(\widetilde\vv)|^{p-2}\nabla\EE(\widetilde\vv)\big)\Vdots
 \nabla\EE(\vv_k{-}\widetilde\vv)
\Big)\,\d \xx\d t
\\[-.1em]&\hspace*{3em}\nonumber
+\int_0^T\!\!\!\int_\varGamma\TRACTION{\cdot}(\vv_k{-}\widetilde\vv)\,\d S\d t
-\int_\varOmega\frac{\varrho_k(T)}2|\vv_k(T)|^2-\frac{\varrho_0}2|\vv_0|^2\,\d \xx
\bigg)%
\\[-.1em]&\nonumber=
 \int_0^T\!\!\!\int_\varOmega
 \Big(\varrho_k
 (\GRAVITY{-}
 \DT\vv_k)
 {\cdot}(\vv_k{-}\widetilde\vv)
  -
\TT_{\varepsilon,k}{:}\nabla(\vv_k{-}\widetilde\vv)
 -\pl_\EE'\zeta(z_k,\EE(\widetilde\vv)){:}\EE(\vv_k{-}\widetilde\vv)
 \\[-.3em]&\hspace*{3em}
-\nu\big(|\nabla\EE(\widetilde\vv)|^{p-2}\nabla\EE(\widetilde\vv)\big)\Vdots
 \nabla\EE(\vv_k{-}\widetilde\vv)
\Big)\,\d \xx\d t
+\int_0^T\!\!\!\int_\varGamma\TRACTION{\cdot}(\vv_k{-}\widetilde\vv)\,\d S\d t\,.
\nonumber\end{align}
Here, we used also $\nabla\vv_{kl}\to\nabla\vv_k$ weakly in
$L^2( I{\times}\varOmega ;\R^{d\times d})$.
 Besides \eqref{conv-of-inirtia.v}, we also used the weak upper semicontinuity
of $
\int_\varOmega
-\frac12\varrho_{kl}(T)|\vv_{kl}(T)|^2\,\d \xx
=\int_\varOmega
-\frac12|\sqrt{\varrho_{kl}(T)}\vv_{kl}(T)|^2\,\d\xx$, together  with the
fact  that, like \eqref{Euler-weak-rho.v(T)}, we have also
\begin{align}
&&&\sqrt{\varrho_{kl}(T)}\vv_{kl}(T)\to\sqrt{\varrho_{k}(T)}\vv_k(T)
&&\text{weakly in $\ L^2(\varOmega;\R^d)$.}&&&&
\label{Euler-weak-rho.v(T)+}\end{align}
This  follows since  $\sqrt{\varrho_{kl}(T)}$ is bounded in
$W^{1,r}(\varOmega)$ and  hinges on the boundedness of the term 
$\pdt{}\sqrt{\varrho_{kl}}=({\rm div}\,\vv_{kl})\sqrt{\varrho_{kl}}
+\vv_{kl}{\cdot}\nabla\varrho_{kl}/\sqrt{\varrho_{kl}}$.

By density arguments, inequality \eq{strong-hyper+} holds for any
$\widetilde\vv\in L^p(I;V_k)$ so that we can substitute
$\widetilde\vv=\vv_k\pm\epsilon\ww$ for $\ww\in L^p(I;V_k)$
with $\ww(T)=0=\ww(0)$. This gives equality in \eq{strong-hyper+} and, dividing
this equality by $\epsilon\ne0$ passing with $\epsilon\to0$,
we obtain the weak formulation of the momentum equation 
\eqref{def-Euler-weak1}, here still  at  its
Galerkin-approximation  level. 
The initial condition $\vv_{kl}(0)=\vv_0$ is kept in the limit, too.

\medskip\noindent{\it Step 5: Further a-priori estimates}. 
   At this point, the only equation which is still discretized is
 the  momentum equation \eq{Euler1-large-CH-reg}.  We
     can perform the ``physical'' test of the  six  equations in \eq{Euler-large-CH-reg} respectively by
 $|\vv_k|^2/2$, $\vv_k$,
$[\widehat\varphi_\varepsilon]_\FF'(\FF_k,z_k)\FF_k^\top$,
$\mu_k$, and $\pdt{}z_k+\nabla\vv_k{\cdot}z_k$,  thus obtaining
estimates   (\ref{est}a--c,e) and \eq{est+} written now for
the weak solution $(\varrho_k,\vv_k,\FF_k,z_k,\mu_k)$ of the
(still semidiscretized) system \eq{Euler-large-CH-reg}.
By comparison, we also obtain an estimate for $\mathcal{N}_k(z_k)=\mu_k-
[\widehat\varphi_\varepsilon]_z'(\FF_k,z_k)$.
Specifically, relying on \eq{Euler-quasistatic-est1-5}
and on the estimates \eq{est-mu} and \eq{est+Fes}, 
we obtain 
\begin{align}\label{est-of-xi}
\big\|\mathcal{N}_k(z_k)\big\|_{L^2(I;H^1(\varOmega))}^{}\le C\,.
\end{align}

\medskip\noindent{\it Step 6: Limit passage  for    
$k\to\infty$}.
We use \TTT sequential weak* compactness and \EEE
the Banach selection principle as in Step~4, now also taking 
\eq{Euler-quasistatic-est3-1+} into account instead of
the estimate in \eq{Euler-quasistatic-est2} which was not uniform in $k$.
 For some not relabeled  subsequence and some $(\varrho,\vv,\FF,z,\mu)$,
we now have
\begin{subequations}\label{Euler-weak++}
\begin{align}
\label{Euler-weak-rho}&\varrho_k\to\varrho&&\hspace*{-11em}
\text{strongly in  $C(I{\times}\barOmega)$}\,
,\\
&\vv_k\to\vv&&\hspace*{-11em}\text{weakly* in $\
L^\infty(I;L^2(\varOmega;\R^d))\cap
L^2(I;W^{2,p}(\varOmega;\R^d))$,}\!\!&&
\label{Euler-weak-v}\\\label{Euler-weak-F}
&\FF_k\to\FF\!\!\!&&\hspace*{-11em}\text{weakly* in $\
L^\infty(I;W^{1,r}(\varOmega;\R^{d\times d}))\,\cap\,
H^1(I;L^2(\varOmega;\R^{d\times d}))$,
}\!\!
\\&&&\hspace*{-7em}\text{and strongly in 
$C(I{\times}\barOmega;\R^{d\times d})$, }
\nonumber\\\label{z-strongly+}
&z_k\to z&&\hspace*{-11em}\text{weakly* in $\ L^\infty(I;L^2(\varOmega))\cap L^2(I;H^1(\varOmega))$}
\\&&&\hspace*{-7em}\text{and strongly in $L^c( I{\times}\varOmega )$ for any $1\le c<2{+}4/d$,}
\nonumber
\\\label{Euler-weak-mu}
&\mu_k\to\mu&&\hspace*{-11em}\text{weakly* in $\ L^\infty(I;L^2(\varOmega))\cap L^2(I;H^1(\varOmega))$}\,,
\\\label{Euler-weak-stress}
&
\TT_{\varepsilon,k}\to\TT_{\varepsilon}=
[\widehat\varphi_\varepsilon]_{\FF}'(\FF,z)\FF^\top\!\!
+\widehat\varphi_\varepsilon(\FF,z)\bbI
\hspace*{-0em}&&\hspace*{0em}\text{strongly in $L^c( I{\times}\varOmega ;\R^{d\times d})$,}
\\&\label{m-strongly+}
\widehat\M(\FF_k,z_k)\to\widehat\M(\FF,z)
&&\hspace*{-6em}\text{strongly in $L^c( I{\times}\varOmega )$ for any $1\le c<\infty$,}
\\&\label{pi-strongly+}
[\widehat\varphi_\varepsilon]_z'(\FF_k,z_k)\to[\widehat\varphi_\varepsilon]_z'(\FF,z)
&&\hspace*{-6em}\text{strongly in $L^c( I{\times}\varOmega )$ for any $1\le c<\infty$.}
\end{align}\end{subequations}

The momentum equation \eq{Euler1-large-CH-reg} (still regularized
by $\varepsilon$)  
is to be treated like in Step~4.  The argument   which led to \eq{rho-v-conv}
is to be now based on the the information about the time derivative
$\pdt{}(\varrho_k\vv_k)$ in a seminorm on $L^{p'}(I;W^{2,p}(\varOmega;\R^d)^*)$
induced by a test by $L^p(I;V_{k_0})$ with $k\ge k_0$, $k_0\in\N$,
or by a Hahn-Banach extension of such time derivatives, cf.\
\cite[Ch.\,8]{Roub13NPDE}. The other terms in
\eq{est-of-DT-rho.v} are bounded in $L^{4/3}(I;L^2(\varOmega;\R^d))$. By a
generalization of the Aubin-Lions compact-embedding theorem,
cf.\ \cite[Lemma 7.7]{Roub13NPDE}, we then obtain 
\begin{align}\label{rho-v-conv+}
&\varrho_{kl}\,\vv_{kl}\to\varrho_{k}\vv_k
&&\hspace*{-1em}\text{strongly in }L^c( I{\times}\varOmega ;\R^d)\ \ \text{ with
any $1\le c<4$,}\,
\end{align}
In fact,  the treatment  of \eq{strong-hyper+}
is to be slightly modified by using first $\widetilde\vv\in H^1(I;V_{k_0})$ and
then, for and $k\ge k_0$, can be used for $(\varrho_k,\vv_k,\TT_{\varepsilon,k},z_k)$
in place of $(\varrho_{kl},\vv_{kl},\TT_{\varepsilon,kl},z_{kl})$ and
$(\varrho,\vv,\TT_{\varepsilon},z)$ in place of
$(\varrho_k,\vv_k,\TT_{\varepsilon,k},z_k)$. Then, by density
arguments, \DELETE{such}  we can resort to some arbitrary  
$\widetilde\vv\in L^p(I;W^{2,p}(\varOmega;\R^d))$.

The limit passage in the semilinear equation
\begin{align}\label{z-k}
\pdt{z_k}+\vv_k{\cdot}\nabla z_k=
{\rm div}\Big(\widehat\M(\FF_k,z_k)\nabla\mu_k
+\big(1{-}\chi_\varepsilon(\FF_k)\big)\nabla z_k\Big)
\end{align}
towards the former equation in \eq{Euler4-large-CH-reg} formulated
weakly is  straightforward  due to (\ref{Euler-weak++}b--e,g).
The limit passage in the equation
\begin{align}\label{mu-k}
\mu_k=[\widehat\varphi_\varepsilon]_z'(\FF_k,z_k)+\mathcal{N}_k(z_k)
\end{align}
towards the variational inequality \eq{def-Euler-weak4} is
simple  by   writing the monotone function $\mathcal{N}_k$ as
the derivative of the Yosida approximation $n_k$
of the indicator function $\delta_{[0,1]}^{}$, i.e.\
$n_k(z)=\min_{0\le\tilde z\le1}|z{-}\tilde z|^2/2$. Thus, using convexity of $n_k$,
\eq{mu-k} can be written as the \UUU variational \EEE inequality
$\int_0^T\int_\varOmega n_k(\widetilde z)+
(\mu_k{-}[\widehat\varphi_\varepsilon]_z'(\FF_k,z_k))
(\widetilde z{-}z_k)\,\d \xx\d t\ge\int_0^T\int_\varOmega n_k(z_k)\,\d \xx\d t$
for $\widetilde z$ valued in $[0,1]$. The limit passage is by 
the convergence (\ref{Euler-weak++}d,e,h) and
the  $\varGamma$-convergence of $n_k$ to $\delta_{[0,1]}^{}$ for $k\to\infty$.

From the calculus  in  \eq{est-of-DT-rho.v}, we can also see the information
$\pdt{}(\varrho\vv)\in L^{p'}(I;W^{2,p}(\varOmega;\R^d)^*)$ while
$\nabla(\varrho\vv)\in L^2(I;L^r(\varOmega;\R^{d\times d}))$  
is like in \eqref{nabla-rho.v}. 

\medskip\noindent{\it Step 7:  Removing the regularization}.
\OK Since \EEE $L^\infty(I;W^{1,r}(\varOmega))\,\cap\,
H^1(I;L^2(\varOmega))$ is embedded in $C( I{\times}\barOmega )$ \OK
for \EEE
$r>d$, 
$\FF$ and its determinant evolve continuously in time, being
valued respectively in $C(\barOmega;\R^{d\times d})$ and $C(\barOmega)$.
Let us recall that, due to \eq{ass-IC} and \OK to \EEE the choice of
$\varepsilon>0$,
the initial condition $\FF_0$ (which \OK is the initial state \EEE
for the $\varepsilon$-regularized system \OK as well) \EEE complies with
the bounds 
in \eq{Euler-quasistatic-est-formal4}. 
Therefore, $\FF$ satisfies these bounds  in \eq{Euler-quasistatic-est-formal4}
not only at $t=0$ but also \OK up to a small positive time.  Indeed,
\EEE  
the $\varepsilon$-regularization of
$1/\det(\cdot)$ and of $\varphi$ is \OK not active, \EEE 
 $(\varrho,\vv,\FF,z,\mu)$ solves the original
nonregularized system \OK for some small time, and  
a-priori bounds \eq{Euler-quasistatic-est-formal4} hold. 
By \OK a \EEE continuation argument, \OK such local-in-time solution
can hence be extended to the whole time interval $I$. In particular,
\EEE 
the $\varepsilon$-regularization \OK remains not active for all
times. \EEE

\medskip\noindent{\it Step 8: Energy balance}.  
\OK Let us conclude by checking that the tests of equations \EEE 
\eq{Euler-large-diff}
respectively by $\vv$, $\Se$, and $\mu$ and of \eq{eq:model1} by
$|\vv|^2$ \OK are legitimate, i.e., rigorously justifiable. These
in turn allow to prove the energy balance
\eq{Euler-swelling-large-energy} integrated over a current time interval
$[0,t]$ via \EEE 
\eq{test-T-nabla-v},
\eq{Euler-swelling-large-dissip}, and \eq{calculus-convective-in-F}. 

\OK The already obtained estimates ensure that \EEE  
$\FF\in L^\infty(I,W^{1,r}(\varOmega;\R^{d\times d}))$, \OK as well as
\EEE 
$\vv\in L^\infty(I;L^2(\varOmega;\R^d))\cap L^2(I;
L^\infty(\varOmega;\R^d))$. \OK From this we deduce \EEE 
$(\nabla\vv)\FF\in L^2(I;L^\infty(\varOmega;\R^{d\times d}))$
and $(\vv{\cdot}\nabla)\FF\in L^2(I;L^r(\varOmega;\R^{d\times d}))$,
\OK and from these two we get  
$\pdt{}\FF=(\nabla\vv)\FF-(\vv{\cdot}\nabla)\FF\in
L^2(I;L^r(\varOmega;\R^{d\times d}))$. 
Thus, the particular terms in \eq{Euler2-diff} are in duality with
$\Se=\widehat\varphi_\FF'(\FF,z)\in L^\infty(I{\times}\varOmega;\R^{d\times d})$.
\OK On the other hand,
we have that \EEE $\varphi(\FF,z)\in L^\infty(I{\times}\varOmega)$ is
in duality \OK with \EEE
${\rm div}\,\vv\in L^4(I{\times}\varOmega))$ \OK and \EEE 
$\mu\in L^2(I; H^1(\varOmega))$ is in duality with
$\DT z\in L^2(I; H^1(\varOmega)^*)$. Thus the tests \eq{test-T-nabla-v} and
\eq{Euler-swelling-large-dissip} \OK can be legitimately performed. \EEE

Similarly, we can see that  $\pdt{}\varrho=-({\rm
  div}\,\vv)\varrho-\vv{\cdot}\nabla\varrho\in L^2(I;L^r(\varOmega))$
is in duality with $|\vv|^2\in L^2(I{\times}\varOmega)$ \OK and \EEE
$\varrho\DT\vv\in L^p(I;W^{2,p}(\varOmega;\R^d)^*)+L^1(I;L^\infty(\varOmega;\R^d))$
is in duality with
$\vv\in L^p(I;W^{2,p}(\varOmega;\R^d))\cap L^\infty(I;L^2(\varOmega;\R^d))$.
\OK Hence, also \EEE the test \eq{calculus-convective-in-F} can be
rigorously \OK performed. \EEE
\end{proof}

\bigskip

{\small

\baselineskip=12pt

\noindent{\it Acknowledgments.}
Support from the \"Osterreichische Austauschdienst-GmbH 
 projects CZ 04/2019 and CZ 01/2021,  from 
M\v SMT \v CR (Ministry of Education of the
Czech Republic) project \linebreak CZ.02.1.01/0.0/0.0/15-003/0000493,
from the institutional support RVO:61388998 (\v CR), and from the
Austrian Science Fund (FWF) projects F\,65, W\,1245, I\,4354, I\,5149,  
is gratefully acknowledged. \UUU The authors express gratitude to the
anonymous referee for valuable comments. \EEE


\begin{thebibliography}{10}

\bibitem{AlCrMa19LRCE}
G.~Alberti, G.~Crippa, and A.~L. Mazzucato.
\newblock Loss of regularity for the continuity equation with non-{L}ipschitz
  velocity field.
\newblock {\em Annals of PDE}, 5:Art.no.9, 2019.

\bibitem{Anan12CHTT}
L.~Anand.
\newblock A {C}ahn-{H}illiard-type theory for species diffusion coupled with
  large elastic-plastic deformations.
\newblock {\em J. Mech. Phys. Solids}, 60:1983--2002, 2012.

\bibitem{BaeSri04DFES}
S.~Baek and A.~R. Srinivasa.
\newblock Diffusion of a fluid through an elastic solid undergoing large
  deformation.
  \newblock {\em Intl. J. Non-Linear Mech.}, 39:201--218, 2004.

\bibitem{Ball}
\UUU J.~M.~Ball.
Global invertibility of Sobolev functions and the interpenetration of matter.
{\it Proc. Roy. Soc. Edinburgh Sect. A}, 88:315--328, 1981.
\EEE

\bibitem{Bana32TOL}
\TTT S.~Banach.
\newblock {\em Th\'eorie des Op\'erations Lin\'eaires}.
\newblock M.~Garasi\'nski, Warszawa, 1932 (Engl. transl. North-Holland, Amsterdam, 1987).\EEE

\bibitem{Biot41GTTS}
M.~A. Biot.
\newblock General theory of three-dimensional consolidation.
\newblock {\em J. Appl. Phys.}, 12:155--164, 1941.

\bibitem{CheAna10CTFP}
S.~A. Chester and L.~Anand.
\newblock A coupled theory of fluid permeation and large deformations for
  elastomeric materials.
\newblock {\em J. Mech. Phys. Solids}, 58:1879--1906, 2010.

\bibitem{CheAna11TMCT}
S.~A. Chester and L.~Anand.
\newblock A thermo-mechanically coupled theory for fluid permeation in
  elastomeric materials: Application to thermally responsive gels.
\newblock {\em J. Mech. Phys. Solids}, 59:1978--2006, 2011.

\bibitem{CuGaTe17SGCF}
M.~Curatolo, S.~Gabriele, and L.~Teresi.
\newblock Swelling and growth: a constitutive framework for active solids.
\newblock {\em Meccanica}, 52:3443--3456, 2017.

\bibitem{Cush97PFHP}
J.~H. Cushman.
\newblock {\em The Physics of Fluids in Hierarchical Porous Media: Angstroms to
  Miles}.
\newblock Springer, Dordrecht, 1997.

\bibitem{Boer05TCMP}
R.~{de Boer}.
\newblock {\em Trends in Continuum Mechanics of Porous Media}.
\newblock Springer, Dordrecht, 2005.

\bibitem{DLReAn14CHTP}
C.V. {Di\,Leo}, E.~Rejovitzky, and L.~Anand.
\newblock A {C}ahn-{H}illiard-type phase-field theory for species diffusion
  coupled with large elastic deformations: Application to phase-separating
  {L}i-ion electrode materials.
\newblock {\em J. Mech. Phys. Solids}, 70:1--29, 2014.

\bibitem{DroChr13CEFE}
A.~D. Drozdov and J.deC. Christiansen.
\newblock Constitutive equations in finite elasticity of swollen elastomers.
\newblock {\em \UUU Internat. \EEE J. Solids \UUU Structures}, \EEE 50:1494--1504, 2013.

\bibitem{DuSoFi10TSMF}
F.~P. Duda, A.~C. Souza, and E.~Fried.
\newblock A theory for species migration in a finitely strained solid with
  application to polymer network swelling.
\newblock {\em J. Mech. Phys. Solids}, 58:515--529, 2010.

\bibitem{FriGur06TBBC}
E.~Fried and M.~E. Gurtin.
\newblock Tractions, balances, and boundary conditions for nonsimple materials
  with application to liquid flow at small-lenght scales.
\newblock {\em Arch. Ration. Mech. Anal.}, 182:513--554, 2006.

\bibitem{GuFrAn10MTC}
M.~E. Gurtin, E.~Fried, and L.~Anand.
\newblock {\em The Mechanics and Thermodynamics of Continua}.
\newblock Cambridge Univ. Press, New York, 2010.

\bibitem{HonWan13PFMS}
W.~Hong and X.~Wang.
\newblock A phase-field model for systems with coupled large deformation and
  mass transport.
  \newblock {\em J. Mech. Phys. Solids}, 61:1281--1294, 2013.

  \UUU
\bibitem{Kroemer}
  S.~Kr\"omer. Global invertibility for orientation-preserving Sobolev
  maps via invertibility on or near the boundary. {\it
    Arch. Ration. Mech. Anal.}, 238:1113--1155, 2020. 
  \EEE

\bibitem{KruRou19MMCM}
M.~Kru\v{z}\'{\i}k and T.~Roub{\'{\i}}{\v{c}}ek.
\newblock {\em Mathematical Methods in Continuum Mechanics of Solids}.
\newblock Springer, Cham/Switzerland, 2019.

\bibitem{LuNaTe13TASI}
A.~Lucantonio, P.~Nardinocchi, and L.~Teresi.
\newblock Transient analysis of swelling-induced large deformations in polymer
  gels.
\newblock {\em J. Mech. Phys. Solids}, 61:205--218, 2013.

\bibitem{Mart19PCM}
Z.~Martinec.
\newblock {\em Principles of Continuum Mechanics}.
\newblock Birkh{\"a}user/Springer, Switzerland, 2019.

\bibitem{Mind64MSLE}
R.~D. Mindlin.
\newblock Micro-structure in linear elasticity.
\newblock {\em \UUU Arch. Rational \EEE Mech. Anal.}, 16:51--78, 1964.

\bibitem{NeNoSi89GSIC}
J.~Ne\v{c}as, A.~Novotn\'y, and M.~\v{S}ilhav\'y.
\newblock Global solution to the ideal compressible heat conductive multipolar
  fluid.
\newblock {\em Comment. Math. Univ. Carolinae}, 30:551--564, 1989.

\bibitem{NeNoSi91GSCI}
J.~Ne\v{c}as, A.~Novotn\'y, and M.~\v{S}ilhav\'y.
\newblock Global solution to the compressible isothermal multipolar fluid.
\newblock {\em J. Math. Anal. Appl.}, 162:223--241, 1991.

\bibitem{NecRuz92GSIV}
J.~Ne\v{c}as and M.~R{\accent23u}\v{z}i\v{c}ka.
\newblock Global solution to the incompressible viscous-multipolar material
  problem.
\newblock {\em J. Elasticity}, 29:175--202, 1992.

\bibitem{Ruzi92MPTM}
M.~R{\accent23u}\v{z}i\v{c}ka.
\newblock Mathematical and physical theory of multipolar viscoelasticity.
\newblock Bonner Mathematische Schriften 233, Bonn, 1992.

\bibitem{Raja07HAMF}
K.~R. Rajagopal.
\newblock On a hierarchy of approximate models for flows of incompressible
  fluids through porous solids.
\newblock {\em Math. Models Meth. Appl. Sci.}, 17:215--252, 2007.

\bibitem{RohLuk17MLDF}
E.~Rohan and V.~Luke\v{s}.
\newblock Modeling large-deforming fluid-saturated porous media using an
  {E}ulerian incremental formulation.
\newblock {\em Adv. Engr. Software}, 113:84--95, 2017.

\bibitem{Roub13NPDE}
T.~Roub{\'\i}{\v{c}}ek.
\newblock {\em Nonlinear Partial Differential Equations with Applications}.
\newblock Birkh\"auser, Basel, 2nd edition, 2013.

\bibitem{Roub21CHEC}
T.~Roub{\'{\i}}{\v{c}}ek.
\newblock Cahn-{H}illiard equation with capillarity in actual deforming
  configurations.
\newblock {\em \UUU Discrete \EEE  Cont. Dynam. Syst. \UUU Ser. \EEE S}, 14:41--55, 2021.

\bibitem{RouSte18TEPR}
T.~Roub{\'{\i}}{\v{c}}ek and U.~Stefanelli.
\newblock Thermodynamics of elastoplastic porous rocks at large strains towards
  earthquake modeling.
\newblock {\em SIAM J. Appl. Math.}, 78:2597--2625, 2018.

\bibitem{RouTom18TMPE}
T.~Roub{\'{\i}}{\v{c}}ek and G.~Tomassetti.
\newblock A thermodynamically consistent model of magneto-elastic materials
  under diffusion at large strains and its analysis.
\newblock {\em Zeit. \UUU Angew. \EEE Math. Phys.}, 69:Art.no.55, 2018.

\bibitem{Stra17MAMP}
B.~Straughan.
\newblock {\em Mathematical Aspects of Multi-Porosity Continua}.
\newblock Springer, Cham/Switzerland, 2017.

\bibitem{Toup62EMCS}
R.~A. Toupin.
\newblock Elastic materials with \UUU couple-stresses. \EEE
\newblock {\em Arch. \UUU Rational \EEE Mech. Anal.}, 11:385--414, 1962.

\bibitem{Silh92MVMS}
M.~\v{S}ilhav\'y.
\newblock Multipolar viscoelastic materials and the symmetry of the coefficient
  of viscosity.
\newblock {\em Appl. Math.}, 37:383--400, 1992.

\end{thebibliography}
\end{document}